\documentclass[aos,preprint]{imsart}

\RequirePackage{amsthm,amsmath,natbib}
\RequirePackage[colorlinks,citecolor=blue,urlcolor=blue]{hyperref}
\RequirePackage{hypernat}

\usepackage{amsfonts}

\usepackage{times}

\usepackage[rgb,x11names]{xcolor}

\arxiv{}
\startlocaldefs

\input DPZ_revision_arXiv.sty

\def\bnote{\bgroup\color{black}}
\def\enote{\egroup}

\endlocaldefs

\begin{document}

\begin{frontmatter}

\title{
Estimation in functional regression for general exponential
families} \runtitle{Functional regression}

\begin{aug}
\author{\fnms{Winston Wei} \snm{Dou}\ead[label=e1]{Winston.Wei.Dou@aya.yale.edu}\thanksref{t1}%
},
\author{\fnms{David} \snm{Pollard}%
\ead[label=e2]{David.Pollard@yale.edu}\thanksref{t2}}
\and
\author{\fnms{Harrison H.} \snm{Zhou}%
\thanksref{t3}%
\ead[label=e3]{Huibin.Zhou@yale.edu}%
\ead[label=u1,url]{http://www.stat.yale.edu/}%
}

\affiliation{Yale University}

\address{%
Statistics Department\\
Yale University\\
\printead{e1}\\
\printead{e2}\\
\printead{e3}\\
\printead{u1}
 }

\thankstext{t1}{Supported in part by NSF FRG grant DMS-0854975}

\thankstext{t2}{Supported in part by NSF grant MSPA-MCS-0528412}

\thankstext{t3}{Supported in part by NSF Career Award DMS-0645676 and NSF FRG grant DMS-0854975}

\end{aug}

\begin{abstract}
This paper studies a class of exponential family models whose
canonical parameters are specified as linear functionals of an
unknown infinite-dimensional slope function. The optimal minimax
rates of convergence for slope function estimation are established.
The estimators that achieve the optimal rates are constructed by
constrained maximum likelihood estimation with parameters whose
dimension grows with sample size. A change-of-measure argument,
inspired by Le~Cam's theory of asymptotic equivalence, is used to
eliminate the bias caused by the nonlinearity of exponential family
models.
\end{abstract}

\begin{keyword}[class=AMS]
\kwd[Primary ]{62J05}
\kwd{60K35}
\kwd[; secondary ]{62G20}
\end{keyword}

\begin{keyword}
\kwd{Approximation of compact operators, Assouad's lemma,
exponential families, functional estimation, minimax rates of
convergence.}
\end{keyword}

\end{frontmatter}

\section{Introduction}  \label{intro}
\bnote There has been extensive exploratory and theoretical study of
{\sl functional data analysis} (FDA) over the past two decades. Two
monographs by \citet{RamsaySilverman2002,RamsaySilverman2005}
provide comprehensive discussions on the methods and applications.

Among many problems involving functional data, slope estimation in
functional linear regression has received substantial attention in
literature: for example, by \citet{CardotFerratySarda2003},
\citet{LiHsing2007}, and \citet{HallHorowitz2007}. In particular,
\citet{HallHorowitz2007} established minimax rates of convergence
and proposed rate-optimal estimators based on spectral truncation
(regression on functional principal components). They showed that
the optimal rates depend on the smoothness of the slope function and
the decay rate of the eigenvalues of the covariance kernel.

In this paper, we study optimal rates of convergence for slope
estimation in functional generalized linear models, for which little
theory is available. We introduce several new technical devices to
overcome the problems caused by nonlinearity of the link function.
To analyze our estimator, we establish a sharp approximation for
maximum likelihood estimators for exponential families parametrized
by linear functions of $m$-dimensional parameters, for an~$m$ that
grows with sample size (see Lemma \ref{AnBn}). We develop a
change-of-measure argument---inspired by ideas from Le~Cam's theory
of asymptotic equivalence of models---to eliminate the effect of
bias terms caused by the nonlinearity of the link function (see
Section \ref{known} and \ref{unknown}).

More precisely, we consider problems where the observed data consist
of independent, identically distributed  pairs $(y_i,\XX_i)$ where
each $\XX_i$ is a Gaussian process indexed by a compact subinterval
of the real line, which with no loss of generality we take to
be~$[0,1]$. We denote the corresponding norm and inner product in
the space $L^2[0,1]$ by~$\norm{\cdot}$ and~$\inner\cdot\cdot$.


We assume, for each $i$, that the random variable $y_i$ conditional
on the process~$\XX_i$, follows a distribution $Q_{\lam_i}$, where
$\{Q_\lam~:~\lam\in \RR\}$ is a one-parameter exponential family. We
take parameter $\lam_i$ to be a linear functional of $\XX_i$ of the
form
\begin{align}
\lam_i &= a +\int_0^1 \XX_i(t)\BB(t)\,dt
\label{lami.def}\\
&\qt{for an unknown constant~$a$ and an unknown $\BB\in L^2[0,1]$.}
\notag
\end{align} \enote
We focus on estimation of~$\BB$ using integrated squared error loss:
$$
L(\BB,\BBhat_n) = \norm{\BB-\BBhat_n}^2 =
\int_0^1\left(\BB(t)-\BBhat_n(t)\right)^2 dt.
$$

Our models are indexed by parameters $f=(K, a,\mu,\BB)$, where $\mu$
is the mean and $K$ is the covariance kernel of the Gaussian
process. 
The universal constant $\alpha$ controls the decay rate of
eigenvalues of kernel $K$ and $\beta$ characterizes the `smoothness'
of the slope function $\BB$. See Definition~\ref{ff.def} (in
Section~\ref{model}) for the precise specification of the parameter
set $\ff = \ff(R,\alpha,\beta)$. The two main results are as
follows.

\Theorem \label{thm.e.upper}{\sl $($Minimax Upper Bound $)$}~ Under
the assumptions stated in Section \ref{model}, there exists an
estimating sequence of~$\BBhat_n$'s for which: for each $\eps>0$
there exists a finite constant~$C_\eps$ such that \beqN
\sup_{f\in\ff}\PP_{n,f}\left\{ \norm{\BB-\BBhat_n}^2>C_\eps
n^{(1-2\beta)/(\alpha+2\beta)}\right\}<\eps \qt{for large enough
$n$.} \eeqN
\endTheorem

\Theorem \label{thm.e.lower}{\sl $($Minimax Lower Bound $)$}~Under
the assumptions stated in Section \ref{model},
 \beqN \liminf_{n\to\infty}n^{(2\beta-1)/(\alpha+2\beta)}\sup_{f\in\ff}\PP_{n,f}
\norm{\BB-\BBhat_n}^2 >0 \qt{for every estimator ~$\{\BBhat_n\}$.}
\eeqN
\endTheorem

Two closely related works in functional data analysis are
\citet{CardotSarda2005} and \citet{MullerStadtmuller2005}, which
provided theory for the functional generalized linear model,
including the rates of convergence for prediction in the random
design case. However, the rate optimalities were not studied. In
addition, \citet{MullerStadtmuller2005} established an upper bound
for rates of convergence assuming the negligibility of the bias due
to the approximation of the infinite-dimensional model by a sequence
of finite-dimensional models, the issue we overcome by using a
change-of-measure argument. In the functional linear regression
setting, \citet{CaiHall2006} and \citet{CrambesKneipSarda2009}
derived optimal rates of convergence for prediction in the fixed and
random design cases. See also, \citet{CardotMasSarda2007} which
derived a CLT for prediction in the fixed and random design cases
and \citet{CardotJohannes2010} which established a minimax optimal
result for prediction at a random design using thresholding
estimators. In a companion study to our paper, \citet[Chapter
5]{Dou2010} considers optimal prediction in functional generalized
linear regressions with an application to the economic problem of
predicting occurrence of recessions from the U.S. Treasury yield
curve.

Our minimax upper bound result (Theorem~\ref{thm.e.upper}) is proved
in Section~\ref{upperproof}. The minimax lower bound result
(Theorem~\ref{thm.e.lower}) is established in
Section~\ref{lowerproof}. The proof of Theorem~\ref{thm.e.upper}
depends on an approximation result (Lemma~\ref{AnBn}) for maximum
likelihood estimators in exponential family models for parameters
whose dimensions change with sample size. As an aid to the reader,
we present our proof of Theorem~\ref{thm.e.upper} in two stages. In
Section~\ref{known}, we assume that both the mean~$\mu$ and the
covariance kernel~$K$ are known. This allows us to emphasize the key
ideas in our proofs without the many technical details that need to
be handled when~$\mu$ and~$K$ are estimated in the natural way. Many
of those details, as summarized in Lemma~\ref{txxn}, involve the
spectral theory of compact operators. We proceed in
Section~\ref{unknown} to the case where~$\mu$ and~$K$ are estimated.
The proofs for the lemmas are collected together in
Section~\ref{lemmaproof}. Some of them invoke the
perturbation-theoretic results collected in the supplemental
Appendix.

\section{Regularity conditions} \label{model}
Let $\{Q_\lam : \lam\in\RR\}$ be a one-parameter exponential family,
\bnote \beq\label{gexpmodel}\id Q_\lam / \id Q_0 = f_{\lam}(y):=
\exp(\lam y -\psi(\lam))\qt{for all}~\lam\in\RR.\eeq Necessarily
$\psi(0)=0$. Remember that $e^{\psi(\lam)} = Q_0 e^{\lam y}$ and
that the distribution $Q_\lam$ has mean $\psidot(\lam)$ and
variance~$\psiddot(\lam)$.
\begin{remark}
We may assume that $\psiddot(\lam)>0$ for every real $\lambda$.
Otherwise we would have $0=\psiddot(\lamb_0)=\var_{\lam_0}(y)=Q_0
f_{\lam_0}(y)(y-\psidot(\lam_0))^2$ for some~$\lam_0$, which would
make $y=\psidot(\lam_0)$ for~$Q_0$ almost all~$y$ and $Q_\lam\equiv
Q_{\lam_0}$ for every~$\lam$.
\end{remark}

We assume:
\begin{enumerate}

\item[(\psiPP)]  \label{psi2}
For each $\eps>0$ there exists a finite constant~$C_\eps$ for
which~$\psiddot(\lam)\le C_\eps\exp(\eps\lam^2)$ for
all~$\lam\in\RR$. Equivalently,  $\psiddot(\lam)\le
\exp\left(o(\lam^2)\right)$ as $|\lam|\to\infty$.

\item[(\psiPPP)] \label{psi3}
There exists an increasing real function $G$ on $\RR^+$ such that
$$
|\psidddot(\lam+h)|\le \psiddot(\lam)G(|h|)\qt{for all $\lam$ and
$h$}.
$$
\Wolog/ we assume $G(0)\ge1$.
\end{enumerate}

As shown in Section~\ref{hellinger}, the assumption (\psiPPP)
implies that \beq h^2(Q_\lam,Q_{\lam+\del})\le
\del^2\psiddot(\lam)\left(1 +|\del|\right)G(|\del|) \qt{for all
$\lam,\del\in\RR$,} \label{hell} \eeq which plays a key role in
analyzing both upper and lower bounds.


We assume the observed data are   iid  pairs $(y_i,\XX_i)$ for $i=1,\dots,n$, where:
\begin{enumerate}

\item[(X)]
Each $\{\XX_i(t):0\le t\le 1\}$ is distributed like~$\{\XX(t):0\le t\le1\}$, a Gaussian process with mean $\mu(t)$ and covariance kernel $K(s,t)$.


\item[(Y)]
$y_i\mid \XX_i\sim Q_{\lam_i}$ with $\lam_i=a+\inner{\XX_i}{\BB}$
for an unknown $\{\BB(t):0\le t\le 1\}$ in $L^2[0,1]$ and $a\in\RR$.
\end{enumerate}

\Definition  \label{ff.def} For real constants $\a>1$ and
$\b>(\a+3)/2$ and~$R>0$, define $\ff= \ff(R,\a,\b)$ as the set of
all $f=(K, a,\mu,\BB)$  that satisfy the following conditions.

\begin{enumerate}


\item[(K)]\label{K}
The covariance kernel is square integrable \wrt/~Lebesgue measure
 and has an eigenfunction expansion (as a compact operator on~$L^2[0,1]$)
$$
K(s,t) = \sumnl_{k\in \NN}\th_k \phi_k(s)\phi_k(t)
$$
where the eigenvalues~$\th_k$ are decreasing with
$
Rk^{-\a}\ge \th_k  \ge \th_{k+1}+ (\a/R) k^{-\a-1}
$.

\item[(a)]\label{a} $|a|\le R$

\item[($\mu$)]
$\norm{\mu}\le R$

\item[($\BB$)]
$\BB$ has an expansion
$
\BB(t) = \sumnl_{k\in\NN}b_k \phi_k(t)
$
with $|b_k|\le Rk^{-\b}$, for the eigenfunctions defined by the kernel~$K$.
\end{enumerate}

\endDefinition

\Remarks The awkward lower bound for $\th_k$ in Assumption~(K)
implies, for all $k<j$, \beq \th_k -\th_j\ge R^{-1}\int_k^j \a
x^{-\a-1}dx  = R^{-1}\left(k^{-\a}-j^{-\a}\right). \label{incr.th}
\eeq If~$K$ and~$\mu$ were known, we would only need the lower
bound~$\th_k\ge R^{-1}k^{-\a}$ and not the lower bound
for~$\th_k-\th_{k+1}$. As explained by
\citet[page~76]{HallHorowitz2007}, the stronger assumption is needed
when one estimates the individual eigenfunctions of~$K$. Note that
the subset~of~$L^2[0,1]$ in which~$\BB$ lies, denoted as $\bb_K$,
depends on~$K$. We regard the need for the stronger assumption on
the eigenvalues and the irksome Assumption~($\BB$) as  artifacts of
the method of proof, but we have not yet succeeded in removing
either assumption.
\endRemarks

More formally, we write $P_{\mu,K}$ for the distribution (a
probability measure on the space~$L^2[0,1]$) of each Gaussian
process~$\XX_i$. The joint distribution of $\XX_1,\dots,\XX_n$ is
then $\PP_{n,\mu,K}=P_{\mu,K}^n$. We identify the $y_i$'s with the
coordinate maps  on~$\RR^n$ equipped with the product measure
$\QQ_{n,a,\BB,\XX_1,\dots,\XX_n}:=\otimes_{i\le n}Q_{\lam_i}$, which
can also be thought of as the conditional joint distribution of
$(y_1,\dots,y_n)$ given $(\XX_1,\dots,\XX_n)$. Thus the~$\PP_{n,f}$
in Theorems~\ref{thm.e.upper} and~\ref{thm.e.lower} can be rewritten
as an iterated expectation,
$$
\PP_{n,f} = \PP_{n,\mu,K}\QQ_{n,a,\BB,\XX_1,\dots,\XX_n} ,
$$
the second expectation on the right-hand side averaging out over
$y_1,\dots,y_n$ for given $\XX_1,\dots,\XX_n$, the first averaging
out over~$\XX_1,\dots,\XX_n$. To simplify notation, we will often
abbreviate $\QQ_{n,a,\BB,\XX_1,\dots,\XX_n}$ to $\QQ_{n,a,\BB}$.


\section{Proof of Theorem \ref{thm.e.upper}}\label{upperproof}
The proof of Theorem \ref{thm.e.upper} will be divided into two
stages. In the first stage, we prove the theorem assuming that the
covariance kernel $K$ is known. This case is relatively simple and
of course artificial, but it captures the essence of the idea of our
proof. In the second stage where $K$ is unknown, we shall show that
using the natural estimate $\widetilde{K}$ as in \cref{tK} will not
affect the result achieved in the first stage. Lemma \ref{txxn} is
to control the gap between the two stages.

In Section \ref{method} we introduce the methodology of constructing
a sequence of estimators achieving the optimal rates of convergence.
In Section~\ref{majorlemma} we state the technical lemmas which
serve as building blocks for establishing the main theorems. Their
proofs are postponed to the Section~\ref{lemmaproof}. In
Section~\ref{known} we prove Theorem~\ref{thm.e.upper} assuming
$\mu$ and $K$ are known, and then in Section~\ref{unknown} we
complete the proof of Theorem~\ref{thm.e.upper} with unknown $\mu$
and $K$.

\subsection{Methodology}\label{method}

Under the assumptions (X) and (K) from Section~\ref{model}, the
process $\XX_i$ admits the eigen decomposition:$$ \XX_i - \mu =\ZZ_i
= \sumnl_{k\in\NN}\zik\phi_k .
$$
The random variables $\zik := \inner{\ZZ_i}{\phi_k}$ are independent
with $\zik\sim N(0,\th_k)$.

Because $\mu$ and~$K$ are unknown, we estimate them in the usual
way: $ \tmu_n(t)=\XXbar(t) = n^{-1}\sumnl_{i\le n} \XX_i(t) $ and
\bAlignL
\tK(s,t) &= (n-1)^{-1}\sumnl_{i\le n} \left(\XX_i(s)-\XXbar(s)\right)\left(\XX_i(t)-\XXbar(t)\right)\label{tK}\\
&= (n-1)^{-1}\sumnl_{i\le n}
\left(\ZZ_i(s)-\ZZbar(s)\right)\left(\ZZ_i(t)- \ZZbar(t)\right),
\notag \eAlignL which has spectral representation
$$
\tK(s,t) = \sumnl_{k\in\NN}\tth_k \tphi_k(s)\tphi_k(t) .
$$
with $\tth_1\geq \tth_2\geq\cdots\geq\tth_{n-1}\geq 0$. In fact we
must have~$\tth_k=0$ for $k\ge n$ because all the
eigenfunctions~$\tphi_k$ corresponding to nonzero~$\tth_k$'s must
lie in the $n-1$-dimensional space spanned by $\{\ZZ_i-\ZZbar:
i=1,2,\dots,n\}$.

Using the first $N$ (as defined in \cref{N}) principal components,
we can approximate the original infinite-dimensional model by the
following sequence of truncated finite-dimensional models: \beqN y_i
| \XX_1,\cdots, \XX_n \sim Q_{\widetilde{\lambda_i}} \eeqN with
$$\widetilde{\lambda}_i = \widetilde{b}_0 + \sumnl_{1\leq j\leq
N}\widetilde{b}_j\widetilde{z}_{i,j},$$ where $\widetilde{b}_0 = a +
\inner{\BB}{\XXbar}$, and $\widetilde{b}_j =
\inner{\BB}{\widetilde{\phi}_j}$ for $j\geq 1$, and
$\widetilde{z}_{i,j} = \inner{\XX_i-\XXbar}{\widetilde{\phi}_j}$. \\

We estimate $\BB$ by \beq\label{bhat}\widehat{\BB} = \sumnl_{j\leq
m}\widehat{b}_j\widetilde{\phi}_j,\eeq where $(\widehat{b}_0,
\cdots, \widehat{b}_N)$ is the conditional MLE for the truncated
model and $m\leq N$. More precisely, $(\widehat{b}_0, \cdots,
\widehat{b}_N)$ is chosen to maximize the following conditional (on
the $\XX_i$'s) log likelihood over $(g_0,g_1,\cdots, g_N)$ in
$\RR^{N+1}$:
\beq\label{bhathat} \ll_n(g_0, g_1,\cdots, g_N) = \sumnl_{i\leq
n}y_i(g_0+\sumnl_{j\leq N}g_j \widetilde{z}_{i,j}) -
\psi(g_0+\sumnl_{j\leq N}g_j \widetilde{z}_{i,j}),\eeq with
\beq\label{m}m\asymp n^{1/(\alpha+2\beta)}\eeq and \beq\label{N}
N\sim n^\zeta\qt{with }(2+2\a)^{-1} > \zeta > (\a+2\b-1)^{-1}. \eeq
Note that $N$ is much larger than~$m$. Such a~$\zeta$ exists because
the assumptions $\a>1$ and $\b>(\a+3)/2$ imply $\a+2\b-1>2+2\a$.

\subsection{Technical lemmas}\label{majorlemma}
We shall first introduce an approximation result for maximum
likelihood estimators in exponential family models for parameters
whose dimensions change with sample size. This lemma combines ideas
from~\citet{Portnoy88} and from~\citet{HjortPollard93}. We write our
results in a notation that makes the applications in
Sections~\ref{known} and~\ref{unknown} more straightforward. The
notational cost is that the parameters are indexed by
$\{0,1,\dots,N\}$. To avoid an excess of parentheses we write $\Np$
for~$N+1$. In the applications $N$ changes with the sample size~$n$
and~$\QQ$ is replaced by~$\QQnaBN$ or~$\tQQnaBN$. For each square
matrix $A$, the spectral norm is defined by $||A||_2:= \sup_{|v|\leq
1}|Av|$ where $|v|$ denotes the $l^2$ norm of vector $v$.

\begin{lemma}\label{AnBn}
Let $Q_\lam$ be the one-parameter exponential family distribution
defined as in \cref{gexpmodel} and satisfying regularity condition
(\psiPPP). Suppose $\xi_1,\dots,\xi_n$ are (nonrandom) vectors
in~$\RRNp$. Suppose $\QQ=\otimes_{i\le n}Q_{\lam_i}$ with
$\lam_i=\xi_i'\gam$ for a
fixed~$\gam=(\gam_0,\gam_1,\dots,\gam_{N})$ in~$\RRNp$. Under~$\QQ$,
the coordinate maps $y_1,\dots,y_n$ are independent random variables
with  $y_i\sim Q_{\lam_i}$.

The log-likelihood for fitting the model is
$$
L_n(g)=\siln(\xi_i'g) y_i -\psi(\xi_i'g)\qt{for $g\in\RRNp$},
$$
which is maximized (over $\RRNp$) at the MLE~$\ghat$ $(=\ghat_n)$.
Suppose~$\xi_i=D\eta_i$ for some nonsingular matrix~$D$, so that
$$
J_n = nDA_nD'\qt{where}~A_n:= \frac1n\sumnl_{i\le
n}\eta_i\eta_i'\psiddot(\lam_i).
$$
If~$B_n$ is another nonsingular matrix for which
\beq\label{assumptionA} \twonorm{A_n-B_n}\le (2\twonorm{\Bni})^{-1}
\eeq and if \beq \maxnl_{i\le n}|\eta_i| \le
\frac{\eps\sqrt{n}/N_+}{ G(1) \sqrt{32\twonorm{ \Bni }  }} \qt{for
some $0<\eps<1$} \label{max.eta} \eeq then for each set of vectors
$\kappa_0,\dots,\kappa_M$ in~$\RRNp$ there is a set
$\yy_{\kappa,\eps}$ with $\QQ \yy_{\kappa,\eps}^c<2\eps$ on which
$$
\sumnl_{0\le j\le M}|\kappa_j'(\ghat-\gam)|^2\le \frac{  6\twonorm{
\Bni } } {n\eps} \sumnl_{0\le j\le M}|D^{-1}\kappa_j|^2.
$$
\end{lemma}

The following approximation result for random matrices will be
invoked in order to apply the Lemma \ref{AnBn} to show Theorem
\ref{thm.e.upper}.
\begin{lemma}
\label{AB} Suppose $\{\eta_{i,k}:i,k\geq 1\}$ are i.i.d. standard
normal random variables. Let \beq\label{A}A_n = n^{-1} \sumnl_{i\le
n}\eta_i\eta_i'\psiddot(\gam'D\eta_i),\eeq where $\gam = (\gamma_0,
\gamma_1, \cdots, \gamma_N)', ~\eta_i =
(1,\eta_{i,1},\dots,\eta_{i,N})' $, and $D = \mbox{diag}(D_0, D_1,
\cdots, D_N)$.
Denote $B_n = \PP A_n$ and assume $\psi$ satisfies condition
(\psiPP). If~$~\sumnl_{k\geq 1}D_k^2\gamma_k^2<\infty$ and
$N=o\left(n^{-1/2}\right)$, it follows that
$\twonorm{\Bni}=O_\ff(1)$ and $\PP \twonorm{A_n-B_n}^2=o_\ff(1)$.

\end{lemma}

The following lemma establishes a bound on the Hellinger distance
between members of an exponential family, which is the key to our
change of measure argument. We write $\hellinger(P,Q)$ for the
Hellinger distance.
\begin{lemma}\label{exp.facts.lemma}
Suppose $\{Q_\lambda~:~\lambda\in \RR\}$ is an exponential family
defined as in \cref{gexpmodel} and satisfies regularity condition
(\psiPPP). Then, $$ \hellinger^2(Q_\lam,Q_{\lam+\del})\le
\del^2\psiddot(\lam)\left(1 +|\del|\right)G(|\del|) ~~~ \forall~~~
\lam,~\del\in \RR.$$
\end{lemma}

The following lemma provides a maximal inequality for
weighted-chi-square variables, which easily leads to maximal
inequalities for Gaussian processes and multivariate normal vectors.
These inequalities will be repeatedly invoked.

\begin{lemma} \label{gaussian} Suppose
$W_i=\sumnl_{k\in\NN}\tau_{i,k}\eta_{i,k}^2$ for $i=1,\dots,n$,
where the $\eta_{i,k}$'s are independent standard normals and the
$\tau_{i,k}$'s are nonnegative constants with $\infty>
T:=\maxnl_{i\le n}\sumnl_{k\in\NN}\tau_{i,k}$. Then
$$
\PP\{\maxnl_{i\le n} W_i > 4T(\log n + x)\} < 2e^{-x} \qt{for each
$x\ge0$.}
$$
\end{lemma}


When we want to indicate that a bound involving constants~$c$, $C$,
$C_1,\dots$  holds uniformly over all models indexed by a set of
parameters~$\ff$, we write $c(\ff)$, $C(\ff)$, $C_1(\ff),\dots$. By
the usual convention for eliminating subscripts, the values of the
constants might change from one paragraph to the next: a
constant~$C_1(\ff)$ in one place needn't be the same as a
constant~$C_1(\ff)$ in another place. For sequences of
constants~$c_n$ that might depend on~$\ff$, we write $c_n=O_\ff(1)$
and $o_\ff(1)$ and so on to show that the asymptotic bounds hold
uniformly over~$\ff$.

\begin{lemma}\label{txxn}
Let $\XX_1, \cdots, \XX_n$ be i.i.d. Gaussian processes satisfying
(X) and (K). Let $m$ and $N$ be integers defined as in \cref{m} and
\cref{N} respectively. Suppose $H_p$ and $\widetilde{H}_p$ are
orthogonal projections operators associated with
$\SPAN\{\phi_1,\cdots, \phi_p\}$ and
$\SPAN\{\widetilde{\phi}_1,\cdots, \widetilde{\phi}_p\}$. Define the
matrix~$\tS:=\diag(\sig_0,\dots,\sig_N)$ with $\sig_0=1$ and
$\sig_k=\sign\big(\inner{\phi_k}{\tphi_k}\big)$ for~$k\ge1$. The key
quantities are:
\begin{enumerate} \itemsep=3pt

\item
$\Del:=\tK-K$

\item $\tD=\diag(1,\tth_1,\dots,\tth_N)^{1/2}$

\item
$\tz_i = (\tz_{i,1},\dots,\tz_{i,N})'$ where
$\tzik=\inner{\ZZ_i}{\tphi_k}$

\item
$\tz_{\cdot }= (\tz_{\cdot 1 },\dots,\tz_{\cdot N })'$ where
$\tzbar{k}= \inner{\ZZbar}{\tphi_k} = n^{-1}\sumnl_{i\le
n}\tz_{i,k}$

\item
$\txi_i = (1,\tz_i'-\tzbar{}')'$ and $\teta_i = D^{-1}\txi_i$. [We
could define $\teta_i = \tD^{-1}\txi_i$ but then we would need to
show that $\tD^{-1}\txi_i\approx D^{-1}\txi_i$. Our definition
merely rearranges the approximation steps.]

\item
$\tgam:=(\tgam_0,\tb_1,\dots,\tb_N)'$ where
$\BB=\sumnl_{k\in\NN}\tb_k\tphi_k$ and $\tgam_0 :=
a+\inner{\BB}{\XXbar}$. [Note that $\lam_i = \tgam_0
+\inner{\BB}{\ZZ_i-\ZZbar}$.]

\item
$\tlamiN =  \tgam_0 + \inner{\tH_N\BB}{\ZZ_i-\ZZbar} =
\txi_i'\tgam$.


\item
$\tA_n = n^{-1}\sumnl_{i\le n}\teta_i\teta_i'\psiddot(\tlam_{i,N})$
\end{enumerate}

For each~$\eps>0$ there exists a set~$\txx_{\eps,n}$, depending
on~$\mu$ and~$K$, with
$$
\sup\NL_\ff\PPnmuK \txx_{\eps,n}^c<\eps \qt{for all large
enough~$n$}
$$
and on which, for some constant~$C_\eps$ that does not depend
on~$\mu$ or~$K$,
\begin{enumerate} \itemsep=3pt

\item \label{txxn.Del}
$\norm{\Del}\le C_\eps n^{-1/2}$

\item \label{txxn.ZZ}
$\maxnl_{i\le n}\norm{\ZZ_i} \le C_\eps\rootlogn$ and
$\norm{\ZZbar}\le C_\eps n^{-1/2}$

\item \label{txxn.Hm}
$\norm{(\tH_m-H_m)\BB}^2 = o_\ff(\rho_n)$

\item \label{txxn.HN}
$\norm{(\tH_N-H_N)\BB}^2 = O_\ff(n^{-1-\nu})$  for some $\nu>0$ that
depends only on~$\a$ and~$\b$

\item \label{txxn.teta}
$ \max_{i\le n}|\teta_i|^2  = o_\ff(\sqrt{n}/N)$

\item \label{txxn.tA}
$\twonorm{\tS\tA_n \tS - A_n} = o_\ff(1)$


\end{enumerate}
\end{lemma}

\subsection{Proof of Theorem \ref{thm.e.upper} with known Gaussian
distribution}\label{known} Initially we suppose that $\mu$ and $K$
are known.
Under~$\QQ_n=\QQ_{n,a,\BB}$, the $y_i$'s are independent, with
$y_i\sim Q_{\lam_i}$ and
$$
\lam_i = a +\inner{\XX_i}{\BB} = b_0 + \sumnl_{k\in\NN}\zik b_k
\qt{where $b_0=a+\inner \mu\BB$.}
$$
Our task is to estimate the $b_k$'s with sufficient accuracy to be
able to estimate $\BB(t) = \sumnl_{k\in\NN}b_k\phi_k(t)$ within an
error of order~$\rho_n = n^{(1-2\b)/(\a+2\b)}$. In fact it will
suffice to estimate the component~$H_m\BB$ of~$\BB$ in the subspace
spanned by $\{\phi_1,\dots,\phi_m\}$ with $m\asymp n^{1/(\a+2\b)}$
because \beq \norm{H_m^\perp \BB}^2 = \sumnl_{k>m}b_k^2 =
O_\ff(m^{1-2\b}) =O_\ff(\rho_n). \label{tail.BB} \eeq

We might try to estimate the coefficients $(b_0,\dots,b_m)$ by
choosing $\ghat=(\ghat_0,\dots,\ghat_m)$ to maximize a conditional
log likelihood over all $g$ in~$\RR^{m+1}$,
$$
\sumnl_{i\le n} y_i\lam_{i,m}-\psi(\lam_{i,m}) \qt{with }\lam_{i,m}
= g_0 + \sumnl_{1\le k\le m}\zik g_k .
$$
To this end we might try to appeal to Lemma~\ref{AnBn} stated at the
beginning of this Section, with $\kappa_j$ equal to the unit vector
with a~$1$ in its~$j$th position for $j\le m$ and $\kappa_j=0$
otherwise. That would give a bound for~$\sumnl_{j\le
m}(\ghat_j-\gam_j)^2$. Unfortunately,  we cannot directly invoke the
Lemma with $N=m$ to estimate $\gam =(b_0,b_1,\dots,b_N)$ when
\begin{align}
\QQ=\QQ_{n,a,\BB} &\AND/ D=\diag(1,\th_1,\dots, \th_N)^{1/2}
\notag\\
\xi_i' = (1,z_{i,1},\dots,z_{i,N}) &\AND/ \eta_i' =
(1,\eta_{i,1},\dots,\eta_{i,N}),\label{xi.gam}
\end{align}
because~$\lam_i\ne \xi_i'\gam$, a bias problem. Note that in this
case $\eta_{i,j} = z_{i,j}/\sqrt{\theta_i}$ for all $i,j$ and hence
the $\eta_{i,j}$'s are i.i.d. standard normal variables.

\Remark We could modify Lemma \ref{AnBn} to allow~$\ell_i=
\xi_i'\gam +\text{bias}_i$, for a suitably small bias term, but at
the cost of extra regularity conditions and a more delicate
argument. The same difficulty arises whenever one investigates the
asymptotics of maximum likelihood with the true distribution outside
the model family.
\endRemark

Instead, we use a two-stage estimation procedure that  eliminates
the bias term by a change of measure conditional on the~$\XX_i$'s.
We shall present the proof in the following three steps.

\subsubsection*{Step 1}From the analysis above, one can see
that the key in our proof is the change-of-measure argument and the
application of Lemma \ref{AnBn}. In this step, we shall construct a
high probability set such that for each realization of the $\XX_i$'s
on the set the assumptions of Lemma \ref{AnBn} are satisfied and the
change-of-measure argument is ready to work.

Define $\xi_i$, $D$, and~$\eta_i$ as in equation~(\ref{xi.gam}).
Then we define matrix $A_n$ as in \cref{A} and choose $B_n :=
\PPnmuK A_n$. Define
$\xx_n=\xx_{\ZZ,n}\cap\xx_{\eta,n}\cap\xx_{A,n}$, where
\begin{align}
\xx_{\ZZ,n} & :=\{\maxnl_{i\le n}\norm{\ZZ_i}^2 \le C_0\log n\}
\label{xxnZZ}\\
\xx_{\eta,n} & :=\{\maxnl_{i\le n}|\eta_i|^2 \le C_0 N\log n\}
\label{xxneta}\\
\xx_{A,n} &:= \{\twonorm{A_n-B_n}\le (2\twonorm{\Bni})^{-1}\}
\label{xxnA}
\end{align}
If we choose a large enough universal constant~$C_0=C_0(\ff)$,
Lemma~\ref{gaussian} ensures that $\PPnmuK\xx_{\ZZ,n}^c\le 2/n$ and
${\PPnmuK\xx_{\eta,n}^c\le 2/n}$ by choosing ${\tau_{i,k} =
\theta_i}$ and ${\tau_{i,k} = \{i\leq N\}}$ respectively for all
$i,k$; and Lemma \ref{AB} shows that
$$
\twonorm{\Bni}=O_\ff(1)\AND/\PPnmuK\twonorm{A_n-B_n}^2=o_\ff(1),
$$
thus $\PPnmuK\xx_{A,n}^c = o_\ff(1)$. And hence,
\beq\label{xn}\PPnmuK\xx_{n}^c\leq \PPnmuK\xx_{\ZZ,n}^c +
\PPnmuK\xx_{\eta,n}^c + \PPnmuK\xx_{A,n}^c = o_\ff(1).\eeq
\subsubsection*{Step 2}  Let us consider the approximate
distribution
$$
\QQnaBN := \otimes_{i\le n}Q_{\lamiN} \qt{with $\lamiN :=
\xi_i'\gam$ and $\gam'=(b_0,b_1,\dots,b_N)$.}
$$
In this step, we show that the divergence caused by
replacing~$\QQnaB$ by $\QQnaBN$ is small enough that it will not
compromise the asymptotic results. In replacing~$\QQnaB$ by
$\QQnaBN$ we eliminate the bias problem but now we have to relate
the probability bounds for~$\QQnaBN$ to bounds involving~$\QQnaB$. A
common control of this divergence is the total variation distance
between $\QQnaBN$ and $\QQnaB$. We shall show that there exists a
sequence of nonnegative constants~$c_n$ of order $o_\ff(\log n)$,
such that \beq \normTV{ \QQ_{n,a,\BB}-\QQnaBN }^2 \le
e^{2c_n}\sumnl_{i\le n}|\lam_i-\lamiN|^2 \qt{on $\xx_n$.}
\label{TV.bound} \eeq To establish inequality~\cref{TV.bound} we use
the bound
$$
\normTV{ \QQ_{n,a,\BB}-\QQnaBN }^2\le
\hellinger^2(\QQ_{n,a,\BB},\QQnaBN) \le \sumnl_{i\le
n}\hellinger^2(Q_{\lam_i},Q_{\lamiN})
$$
By Lemma~\ref{exp.facts.lemma}
$$
\hellinger^2(Q_{\lam_i},Q_{\lamiN}) \le
\del_i^2\psiddot(\lam_i)\left(1+|\del_i|\right)g(|\del_i|)
$$
where \bAlign
|\del_i| =|\lam_i-\lamiN| &= |\inner{\ZZ_i}{\BB} -\inner{H_N\ZZ_i}{\BB}|\\
&=|\inner{\ZZ_i}{H_N^\perp\BB}| \le \norm{\ZZ_i}\norm{H_N^\perp\BB}\\
&\le  O_\ff\left(\sqrt{N^{1-2\b}\log n }\right) = o_\ff(1) \eAlign
Thus all the $\left(1+|\del_i|\right)g(|\del_i|)$ factors can be
bounded by a single~$O_\ff(1)$ term.

For $(a,\BB,\mu,K)\in\ff(R,\a,\b)$ and with the $\norm{\ZZ_i}$'s
controlled by~$\xx_n$,
$$
|\lam_i| \le |a| +(\norm{\mu}+\norm{\ZZ_i})\norm{\BB} \le
C_2\sqrt{\log n}
$$
for some constant~$C_2=C_2(\ff)$. Assumption~(\psiPP) then ensures
that all the  $\psiddot(\lam_i)$  are bounded by a single
$\exp\left(o_\ff(\log n)\right)$ term.

\subsubsection*{Step 3}
 On the set $\xx_n$, we can apply Lemma~\ref{AnBn} directly with $\QQ = \QQnaBN$, because inequality
\cref{assumptionA} holds by construction and inequality
\cref{max.eta} holds for large enough~$n$ because
$$
\max\NL_{i\le n}|\eta_i|^2 \le O_\ff(N\log n) = o_\ff(\sqrt{n}/N) .
$$

Estimate $\gam$ by the $\ghat=(\ghat_0,\dots,\ghat_N)$ defined in
Lemma~\ref{AnBn}. Thus, the estimator in Theorem \ref{thm.e.upper}
is $\BBhat_n = \sumnl_{1\le k\le m}\ghat_k\phi_k$. For each
realization of the~$\XX_i$'s in~$\xx_n$, Lemma~\ref{AnBn} gives a
set $\yy_{m,\eps}$ with $\QQnaBN\yy_{m,\eps}^c<2\eps$ on which
$$
\sumnl_{1\le k\le m}|\ghat_k-\gam_k|^2 =
O_\ff\left(n^{-1}\sumnl_{1\le k\le m}\th_k^{-1}\right) =
O_\ff(m^{1+\a}/n)=O_\ff(\rho_n),
$$
which implies
$$
\norm{\BBhat_n-\BB}^2 = \sumnl_{1\le k\le m}|\ghat_k-\gam_k|^2 +
\sumnl_{k>m}b_k^2
 =  O_\ff(\rho_n).
$$

From the inequality \cref{TV.bound} it follows, for a large enough
constant~$C_\eps$, that \bAlign
\PPnmuK&\QQnaB\{\norm{\BBhat_n-\BB}^2 > C_\eps\rho_n \}  \\
&\le \PP_{n,\mu,K}\xx_n^c +\PP_{n,\mu,K}\xx_n\left(\normTV{\QQnaB -\QQnaBN} + \QQnaBN \yy_{m,\eps}^c \right)\\
&\le o_\ff(1) +2\eps + e^{c_n}\left(\sumnl_{i\le
n}\PP_{n,\mu,K}|\lam_i-\lam_{i,N}|^2\right)^{1/2} . \eAlign By
construction,
$$
\lam_i-\lam_{i,N} = \sumnl_{k> N}\zik b_k
$$
with the $\zik$'s independent and $\zik\sim N(0,\th_k)$. Thus
$$
\sumnl_{i\le n}\PP_{n,\mu,K}|\lam_i-\lamiN|^2 \le n\sumnl_{k> N}
\th_kb_k^2 = O_\ff(nN^{1-\a-2\b})= o_\ff(e^{-2c_n})
$$
because $\zeta>(\a+2\b-1)^{-1}$. That is, we have an estimator that
achieves the $O_\ff(\rho_n)$ minimax rate.

\subsection{Proof of Theorem \ref{thm.e.upper} with unknown Gaussian
distribution}\label{unknown} As before, most of the analysis will be
conditional on the~$\XX_i$'s lying in a set with high probability on
which the various estimators and other random quantities are well
behaved. Remember $\txx_{\eps,n}$ is the high probability set
defined in Lemma~\ref{txxn}. For the key quantities defined in
Lemma~\ref{txxn}, we shall keep their notations unchanged in this
section for the purpose of making the application more
straightforward.

As before, the component of~$\BB$ orthogonal to
$\SPAN\{\tphi_1,\dots,\tphi_m\}$ causes no trouble because
$$
\norm{\BBhat-\BB}^2 = \sumnl_{1\leq k\leq m}(\ghat_k-\tgam_k)^2
+\norm{\tH_m^\perp\BB}^2
$$
and, by~Lemma~\ref{txxn} part \ref{txxn.Hm},
$$
\norm{\tH_m^\perp\BB}^2 \le 2\norm{H_m^\perp\BB}^2
+2\norm{(\tH_m-H_m)\BB}^2 = O_\ff(\rho_n) \qt{on $\txx_{\eps,n}$}.
$$
To handle 
$\sumnl_{1\leq k\leq m}(\ghat_k-\tgam_k)^2$, invoke Lemma~\ref{AnBn}
for~$\XX_i$'s in~$\txx_{\eps,n}$, with $\eta_i$ replaced
by~$\teta_i$ and $A_n$ replaced by~$\tA_n$ and $B_n$ replaced
by~$\tB_n = \tS B_n\tS$, the same $B_n$ and~$D$ as before, and $\QQ$
equal to
$$
\tQQnaBN = \otimes_{i\le n} Q_{\tlamiN}  .
$$
to get a set $\tyy_{m,\eps}$ with $\tQQnaBN\tyy_{m,\eps}^c<2\eps$ on
which 
$\sumnl_{1\leq k\leq m}(\ghat_k-\tgam_k)^2$.  The conditions of
Lemma~\ref{AnBn} are satisfied on~$\txx_{\eps,n}$, because
of~Lemma~\ref{txxn} part \ref{txxn.teta} and
$$
\twonorm{\tA_n-\tB_n} \le  \twonorm{\tA_n- \tS A_n \tS}+\twonorm{\tS
A_n\tS-\tS B_n\tS} =o_\ff(1).
$$

To complete the proof it suffices to show that
$\normTV{\QQnaBN-\tQQnaBN}$ tends to zero. First note that \bAlign
\tlam_{i,N}-\lam_{i,N} &= a +\inner{\BB}{\XXbar}
+\inner{\tH_N\BB}{\ZZ_i-\ZZbar} -a-\inner{\BB}{\mu} -
\inner{H_N\BB}{\ZZ_i}
\\
&= \inner{\tH_N^\perp \BB}{\ZZbar} -\inner{H_N^\perp
\BB}{\ZZbar}+\inner{H_N^\perp \BB}{\ZZbar} +
\inner{\tH_N\BB-H_N\BB}{\ZZ_i} \eAlign which implies that,
on~$\txx_{\eps,n}$, \bAlignL |\tlam_{i,N}-\lam_{i,N}|^2 &\le
2|\inner{H_N^\perp \BB}{\ZZbar}|^2 +
2\norm{\tH_N\BB-H_N\BB}^2\left(\norm{\ZZ_i} + \norm{\ZZbar}\right)^2
\notag\\
&\le O_\ff(N^{1-2\b})C_\eps^2 n^{-1} +
O_\ff(n^{-1-\nu})C_\eps^2\left(n^{-1/2}+\rootlogn\right)^2
\notag\\
&= O_\ff(n^{-1-\nu'})\qt{for some $0<\nu'<\nu$.} \label{tlam.lam}
\eAlignL Now argue as in the step 2 of the proof for the case of
known $K$:~ on~$\txx_{\eps,n}$, \bAlign \normTV{\tQQnaBN-\QQnaBN}^2
&\le \sumnl_{i\le n}h^2\left(Q_{\tlamiN},Q_{\lamiN}\right)
\\
&\le \exp\left(o_\ff(\log n)\right)\sumnl_{i\le n}|\tlamiN-\lamiN|^2
=o_\ff(1) . \eAlign Finish the argument as before, by splitting into
contributions from~$\txx_{\epsilon,n}^c$ and
${\txx_{\epsilon,n}\cap\tyy_{m,\eps}^c}$ and
${\txx_{\epsilon,n}\cap\tyy_{m,\eps}}$.

\section{Proof of Theorem \ref{thm.e.lower}}\label{lowerproof}
We apply a slight variation on Assouad's Lemma---combining ideas
from \citet{Yu97} and from \citet[Section~24.3]{vanderVaart98}---to
establish the minimax lower bound result in
Theorem~\ref{thm.e.lower}.

We consider behavior only for $\mu=0$ and $a=0$, for a fixed~$K$
with spectral decomposition $\sumnl_{j\in\NN}\th_j
\phi_j\otimes\phi_j$. For simplicity we abbreviate $\PP_{n,0,K}$
to~$\PP$. Let $J=\{m+1,m+2,\dots,2m\}$ and $\Gam=\{0,1\}^J$.
Let~$\b_j=R j^{-\b}$. For each $\gam$ in~$\Gam$ define $\BB_\gam=
\eps\sumnl_{j\in J}\gam_j\b_j \phi_j$, for a small~$\eps>0$ to be
specified, and write $\QQ_\gam$ for the product measure
$\otimes_{i\le n}Q_{\lam_i(\gam)}$ with
$$
\lam_i(\gam)=\inner{\BB_\gam}{\ZZ_i} = \eps\sumnl_{j\in J}\gam_j\b_j
z_{i,j}.
$$

For each~$j$ let $\Gam_j=\{\gam\in\Gam: \gam_j=1\}$ and let $\psi_j$
be the bijection on~$\Gam$ that flips the $j$th coordinate but
leaves all other coordinates unchanged. Let $\pi$ be the uniform
distribution on~$\Gam$, that is, $\pi_\gam=2^{-m}$ for each~$\gam$.

For each estimator $\BBhat=\sumnl_{j\in\NN}\bhat_j\phi_j$ we have $
\norm{\BB_\gam-\BBhat}^2 \ge \sumnl_{j\in
J}\left(\gam_j\b_j-\bhat_j\right)^2 $ and so \bAlignL
\sup_{\ff}\PP_{n,f}\norm{\BB-\BBhat}^2 &\ge \sumnl_{\gam\in\Gam}
\pi_\gam\sumnl_{j\in J}\PP
\QQ_\gam\left(\eps\gam_j\b_j-\bhat_j\right)^2
\notag\\
&= 2^{-m}\sumnl_{j\in J}\sumnl_{\gam\in\Gam_j}\PP \left(\QQ_\gam
(\eps\b_j-\bhat_j)^2 +\QQ_{\psi_j(\gam)}(0-\bhat_j)^2\right)
\notag\\
&\ge 2^{-m}\sumnl_{j\in J}\sumnl_{\gam\in\Gam_j}\tfrac14(\eps\b_j)^2
\PP\norm{\QQ_\gam \wedge \QQ_{\psi_j(\gam)}}, \label{assouad}
\eAlignL the last lower bound coming from the fact that
$$
(\eps\b_j-\bhat_j)^2 + (0-\bhat_j)^2\ge \tfrac14(\eps\b_j)^2 \qt{for
all $\bhat_j$}.
$$
We assert that, if $\eps$ is chosen appropriately, \beq
\min\NL_{j,\gam} \PP\norm{\QQ_\gam \wedge \QQ_{\psi_j(\gam)}} \text{
stays bounded away from zero as $n\to\infty$}, \label{affinity} \eeq
which will ensure that the lower bound in~\cref{assouad} is
eventually larger than a constant multiple of $ \sumnl_{j\in
J}\b_j^2 \ge c\rho_n $ for some constant $c>0$. The inequality in
Theorem~\ref{thm.e.lower} will then follow.

To prove~\cref{affinity}, consider a~$\gam$ in~$\Gam$ and the
corresponding~$\gam'=\psi_j(\gam)$. By virtue of the inequality
$$
\norm{\QQ_\gam\wedge \QQ_{\gam'}} = 1-\normTV{\QQ_\gam- \QQ_{\gam'}}
\ge 1 - \left(2\wedge\sumnl_{i\le
n}h^2(Q_{\lam_i(\gam)},Q_{\lam_i(\gam')})\right)^{1/2}
$$
it is enough to show that \beq
\limsup\NL_{n\to\infty}\max\NL_{j,\gam}\PP\left(2\wedge\sumnl_{i\le
n}h^2(Q_{\lam_i(\gam)},Q_{\lam_i(\gam')})\right) <1 . \label{pph2}
\eeq Define $\xx_n = \{\max\NL_{i\le n}\norm{\ZZ_i}^2 \le C_0\log
n\}$. Based on Lemma \ref{gaussian}, we know that $\PP\xx_n^c=o(1)$
with the constant~$C_0$ large enough. On~$\xx_n$ we have
$$
|\lam_i(\gam)|^2 \le \sumnl_{j\in J}\b_j^2 \norm{Z_i}^2 =
O(\rho_n)\log n =o(1)
$$
and, by inequality~\cref{hell},
$$
h^2(Q_{\lam_i(\gam)},Q_{\lam_i(\gam')}) \le O_\ff(1)
|\lam_i(\gam)-\lam_i(\gam')|^2 \le \eps^2 O_\ff(1)\b_j^2 z_{i,j}^2 .
$$
We deduce that \bAlign \PP\left( 2\wedge\sumnl_{i\le n}
h^2(Q_{\lam_i(\gam)},Q_{\lam_i(\gam')})\right) &\le
2\PP\xx_n^c+\sumnl_{i\le n}\eps^2 O_\ff(1)\b_j^2\PP \xx_nz_{i,j}^2
\\
&\le o(1)+ \eps^2 O(1) n\b_j^2\th_j . \eAlign The choice of~$J$
makes $\b_j^2\th_j\le R^2m^{-\a-2\b}\sim R^2/n$.
Assertion~\cref{pph2} follows.

\section{Proof of technical lemmas}\label{lemmaproof}
\subsection{Proof of Lemma \ref{AnBn}} We need to first show the
following lemma. Define
\begin{enumerate}
\item
$w_i :=J_n^{-1/2}\xi_i$, an element of~$\RRNp$
\item
$W_n= \siln w_i\left(y_i-\psidot(\lam_i)\right)$, an element
of~$\RRNp$
\end{enumerate}

Notice that $\QQ W_n=0$ and $ \var_\QQ(W_n)=\siln
w_iw_i'\psiddot(\lam_i) = I_\Np$ and
$$
\QQ|W_n|^2 =\trace \left(\var_\QQ(W_n)\right) =\Np .
$$

\Lemma  \label{mle.approx} Suppose $0<\eps_1 \le 1/2$ and
$0<\eps_2<1$  and
$$
\maxnl_{i\le n}|w_i| \le \frac{\eps_1\eps_2}{2G(1)\Np} \qt{with $G$
as in Assumption~(\psiPPP).}
$$
Then $\ghat=\gam+J_n^{-1/2}(W_n+r_n)$ with $|r_n|\le \eps_1$ on the
set $\{|W_n|\le \sqrt{N_+/\eps_2}\}$, which has $\QQ$-probability
greater than~$1-\eps_2$.
\endLemma
\begin{proof}
 The equality $\QQ|W_n|^2=\Np$ and Tchebychev give
$$ {\QQ\{|W_n|> \sqrt{\Np/\eps_2}\}\le\eps_2}. $$
Reparametrize by defining $t=J_n^{1/2}(g-\gam)$. The concave
function
$$
\ll_n(t) :=L_n(\gam+J_n^{-1/2}t)-L_n(\gam) =\sumnl_{i\le n}y_iw_i't
+\psi(\lam_i)-\psi(\lam_i+w_i't)
$$
is maximized at $\that=J_n^{1/2}(\ghat-\gam)$. It  has derivative
$$
\lldot_n(t) = \siln w_i\left(y_i -\psidot(\lam_i+w_i't)\right) .
$$
For a fixed unit vector~$u\in\RRNp$ and a fixed~$t\in\RRNp$,
consider the real-valued function of the real variable~$s$,
$$
H(s) := u'\lldot_n(st) = \siln u'w_i\left(y_i
-\psidot(\lam_i+sw_i't)\right),
$$
which has derivatives \bAlign
\Hdot(s) &= -\siln (u'w_i) (w_i't)\psiddot(\lam_i+sw_i't)\\
\Hddot(s) &= -\siln (u'w_i) (w_i't)^2\psidddot(\lam_i+sw_i't)  .
\eAlign Notice that $H(0)=u'W_n$ and $\Hdot(0)= -u'\siln w_i
w_i'\psiddot(\lam_i)t=-u't$.

Write $M_n$ for $\maxnl_{i\le n}|w_i|$. By virtue of
Assumption~(\psiPPP), \bAlign
|\Hddot(s)|&\le \siln|u'w_i| (w_i't)^2 \psiddot(\lam_i)G\left(|sw_i't|\right)\\
&\le M_n G\left(M_n|st|\right)t'\siln w_iw_i'\psiddot(\lam_i)t\\
&= M_n G\left(M_n|st|\right)|t|^2 . \eAlign By Taylor expansion, for
some $0<s^*<1$,
$$
|H(1)- H(0)-\Hdot(0)| \le \tfrac12 |\Hddot(s^*)| \le \tfrac12  M_n
G\left(M_n|t|\right)|t|^2 .
$$
That is, \beq \left|u'\left(\strut\lldot_n(t)-W_n+t\right)\right|
\le \tfrac12 M_n G\left(M_n|t|\right)|t|^2. \label{u'll} \eeq
Approximation~\cref{u'll} will  control  the behavior of  $\tll(s)
:=\ll_n(W_n+su)$, a concave function of the real argument~$s$, for
each unit vector~$u$. By concavity, the derivative
$$
\tlldot(s) = u'\lldot_n(W_n+su) = -s +R(s)
$$
is a decreasing function of~$s$ with
$$
|R(s)|\le \tfrac12 M_n G\left(M_n|W_n+su|\right)|W_n+su|^2
$$
On the set $\{|W_n|\le \sqrt{\Np/\eps_2}\}$ we have
$$
|W_n\pm\eps_1 u|\le \sqrt{\Np/\eps_2} +\eps_1.
$$
Thus
$$
M_n |W_n\pm\eps_1 u|\le
\frac{\eps_1\eps_2}{2G(1)\Np}\left(\sqrt{\Np/\eps_2}
+\eps_1\right)<1,
$$
implying \bAlign |R(\pm\eps_1)| &\le \tfrac12 M_n G(1)|W_n\pm \eps_1
u|^2
\\
&\le \frac{\eps_1\eps_2}{G(1)\Np}\left(\Np/\eps_2 +\eps_1^2\right)
\\
&\le \eps_1\left(1+\eps_1^2\eps_2/\Np\right) <\tfrac58\eps_1 .
\eAlign Deduce that \bAlign
\tlldot(\eps_1) &= -\eps_1+R(\eps_1)\le -\tfrac38\eps_1\\
\tlldot(-\eps_1) &= \eps_1+R(-\eps_1)\ge \tfrac38\eps_1 \eAlign The
concave function $s\mapsto\ll_n(W_n+su)$ must achieve its maximum
for some $s$ in the interval~$[-\eps_1,\eps_1]$, for each unit
vector~$u$. It follows that $|\that-W_n|\le \eps_1$.
\end{proof}
First we establish a bound on the spectral distance between~$\Ani$
and~$\Bni$. Define $H=\Bni A_n-I$.
Then$\twonorm{H}\le\twonorm{\Bni}\twonorm{A_n-B_n}\le 1/2$, which
justifies the expansion
$$
\twonorm{\Ani-\Bni} = \twonorm{\left((I+H)^{-1}-I\right)\Bni} \le
\sumnl_{j\ge 1}\twonorm{H}^k\twonorm{\Bni} \le \twonorm{\Bni} .
$$
As a consequence, $\twonorm{ \Ani } \le 2\twonorm{ \Bni }$.

Choose $\eps_1=1/2$ and $\eps_2=\eps$ in Lemma~\ref{mle.approx}. The
bound on $\maxnl_{i\le n}|\eta_i|$ gives the bound on $\maxnl_{i\le
n}|w_i|$ needed by the Lemma:
$$
n|w_i|^2 = \eta_i'D(J_n/n)^{-1}D\eta_i = \eta_i'\Ani\eta_i \le
\twonorm{ \Ani }|\eta_i|^2 .
$$

Define~$K_j :=J_n^{-1/2}\kappa_j$, so that $
|\kappa_j'(\ghat-\gam)|^2 \le 2(K_j'W_n)^2 + 2(K_j'r_n)^2 $. By
Cauchy-Schwarz,
$$
\sumnl_{j} (K_j'r_n)^2 \le \sumnl_{j}|K_j|^2 |r_n|^2 = U_\kappa
|r_n|^2
$$
where \bAlign U_\kappa:=\sumnl_{j}\kappa_j'J_n^{-1}\kappa_j &=
\sumnl_{j}n^{-1}(D^{-1}\kappa_j)'\Ani D^{-1}\kappa_j
\\
&\le  2n^{-1}\twonorm{\Bni} \sumnl_{j}|D^{-1}\kappa_j|^2 . \eAlign
For the contribution $V_\kappa := \sumnl_{j}|K_j'W_n|^2$ the
Cauchy-Schwarz bound is too crude. Instead, notice that $\QQ
V_\kappa =U_\kappa$, which ensures that the complement of the set
$$
\yy_{\kappa,\eps} := \{|W_n|\le \sqrt{\Np/\eps}\} \cap \{V_\kappa
\le U_\kappa/\eps\}
$$
has $\QQ$ probability less that~$2\eps$.  On the
set~$\yy_{\kappa,\eps}$,
$$
\sumnl_{0\le j\le N}|\kappa_j'(\ghat-\gam)|^2 \le 2V_\kappa
+2U_\kappa|r_n|^2 \le 3U_\kappa/\eps .
$$
The asserted bound  follows.

\subsection{Proof of Lemma \ref{AB}}
Throughout this subsection abbreviate $\PP_{n,\mu,K}$ to~$\PP$.

The matrix $A_n$ is an average of~$n$ independent random matrices
each of which is distributed like~$\nn\nn'\psiddot(\gam'D\nn)$,
where $\nn'=(\nn_0,\nn_1,\dots,\nn_N)$ with $\nn_0\equiv1$ and the
other $\nn_j$'s are independent~$N(0,1)$'s.  Moreover, by rotational
invariance of the spherical normal, we may assume with no loss of
generality that $\gam'D\nn = \abar+\kappa\nn_1$, where
$$
\kappa^2 = \sumnl_{k=1}^N D_k^2 b_k^2  = O_\ff(1).
$$
Thus
$$
B_n =\PP \nn\nn'\psiddot(\abar+\kappa\nn_1) = \diag(F, r_0 I_{N-1})
$$
where
$$
r_j := \PP \nn_1^j\psiddot(\abar+\kappa\nn_1) \AND/ F =
\begin{bmatrix}r_0&r_1\\ r_1&r_2\end{bmatrix}.
$$
The block diagonal form of~$B_n$ simplifies calculation of spectral
norms. \bAlign
\twonorm{\Bni} &= \twonorm{\diag(F^{-1},r_0^{-1}I_{N-1})} \\
&\le \max\left(\twonorm{F^{-1} }, \twonorm{r_0^{-1}I_{N-1}}\right)
\le \max\left(\frac{r_0+r_2}{r_0r_2-r_1^2}, r_0^{-1}\right) .
\eAlign Assumption~(\psiPP)  ensures that both $r_0$ and $r_2$ are
$O_\ff(1)$.

Continuity and strict positivity of~$\psiddot$, together with
$\max(|\abar|,\kappa)=O_\ff(1)$, ensure that $ c_0 :=
\inf\NL_{\abar,\kappa}\inf\NL_{|x|\le 1}\psiddot(\abar+\kappa x)>0
$. Thus
$$
\sqrt{2\pi} r_0 \ge c_0\int_{-1}^{+1} e^{-x^2/2}dx >0
$$
Similarly \bAlign \sqrt{2\pi}(r_0r_2-r_1^2)
&= \sqrt{2\pi}r_0\PP\psiddot(\abar+\kappa\nn_1)(\nn_1-r_1/r_0)^2\\
&\ge  c_0r_0\int_{-1}^{+1} (x-r_1/r_0)^2e^{-x^2/2}dx \ge
c_0r_0\int_{-1}^{+1} x^2e^{-x^2/2}dx . \eAlign It follows that
$\twonorm{\Bni}=O_\ff(1)$.

The random matrix $A_n-B_n$ is an average of~$n$ independent random
matrices each distributed like~$\nn\nn'\psiddot(\abar+\kappa\nn_1)$
minus its expected value. Thus
$$
\PP\twonorm{A_n-B_n}^2\le \PP\Fnorm{A_n-B_n}^2 = n^{-1}\sumnl_{0\le
j,k \le N}\var\left(\nn_j\nn_k \psiddot(\abar+\kappa\nn_1)\right).
$$
Assumption~(\psiPP)  ensures that each summand is~$O_\ff(1)$, which
leaves us with a $O_\ff(N^2/n) =o_\ff(1)$ upper bound.
\subsection{Proof of Lemma \ref{exp.facts.lemma}}\label{hellinger}
Let us temporarily write $\lam'$ for $\lam+\del$ and write $\barlam$
for ${(\lam+\lam')/2=\lam+\del/2}$. \bAlign
1-\tfrac12\hellinger^2(Q_\lam,Q_{\lam'}) &=\int\sqrt{f_\lam(y)f_{\lam'}(y)}\\
&= \int \exp\left(\barlam y -\tfrac12\psi(\lam)-\tfrac12\psi(\lam')\right)\\
&= \exp\left(\psi(\barlam) -\tfrac12\psi(\lam)-\tfrac12\psi(\lam')\right)\\
&\ge 1 + \psi(\barlam) - \tfrac12\psi(\lam) -\tfrac12\psi(\lam')
\eAlign That is,
$$
\hellinger^2(Q_\lam,Q_{\lam'})\le  \psi(\lam)+\psi(\lam+\del) -
2\psi(\lam+\del/2).
$$
By Taylor expansion in~$\del$ around~$0$, the right-hand side is
less than \bAlign
 \tfrac14\del^2\psiddot(\lam) &+\tfrac16\del^3\left(\psidddot(\lam+\del^*)-\tfrac18\psidddot(\lam-\del^*/2)\right)
\eAlign where $0<|\del^*|<|\del|$.  Invoke inequality~\psiPPP twice
to bound the coefficient of $\del^3/6$ in absolute value by
$$
\psiddot(\lam)\left(\strut G(|\del|) +\tfrac18 G(|\del|/2)\right)
\le \tfrac98\psiddot(\lam)G(|\del|) .
$$
The stated bound simplifies some unimportant constants.

\subsection{Proof of Lemma \ref{gaussian}}
\Wolog/, let us suppose $T=1$. For $s=1/4$, note that \bAlign \PP
\exp(sW_i) &=\prodnl_{k\in\NN}(1-2s\tau_{i,k})^{-1/2} \le
\exp\left(\sumnl_{k\in\NN}s\tau_{i,k}\right)\le e^{1/4} \eAlign by
virtue of  the inequality $-\log(1-t)\le 2t$ for $|t|\le 1/2$. With
the same~$s$, it then follows that \bAlign
\PP\{\maxnl_{i\le n}& W_i > 4(\log n + x)\}\\
&\le \exp\left(-4s(\log n + x)\right)\PP\exp\left(\maxnl_{i\le n} sW_i \right)\\
&\le e^{-x}\frac1n\sumnl_{i\le n}\PP \exp(sW_i). \eAlign The $2$ is
just a  clean upper bound for~$e^{1/4}$.

\subsection{Proof of Lemma \ref{txxn}} We shall first show some
preliminary results that will be used in the main proof throughout
Sections \ref{txxnDel} to \ref{txxntA}. In this section, for
notational simplicity, we write $\sumnl_j^*$ for $\sumnl_{j\neq k}$.

Many of the inequalities in this section involve sums of functions
of the $\th_j$'s. The following result will save us a lot of
repetition. To simplify the notation, we drop the subscripts
from~$\PPnmuK$.

\Lemma  \label{weights} \

\begin{enumerate}
\ritem(i) For each $r\ge 1$ there is a constant $C_r=C_r(\ff)$ for
which
$$
\kappa_k(r,\gam) := \sumnl_{j\in\NN}\{j\ne
k\}\frac{j^{-\gam}}{|\th_j-\th_k|^r} \le
\begin{cases}
C_r \left(1+k^{r(1+\a)-\gam}\right)& \text{if $r>1$}\\
C_1 \left(1+k^{1+\a-\gam}\log k\right)&\text{if $r=1$}
\end{cases}
$$
\ritem(ii) For each $p$,
$$
\sumnl_{k\le p}\sumnl_{j>p} \frac{ k^{-\a-2\b}j^{-\a} }{
|\th_k-\th_j|^2} = O_\ff(p^{1-\a})
$$

\end{enumerate}

\endLemma
\begin{proof} For (i), argue in the same way as
\citet[page~85]{HallHorowitz2007}, using the lower bounds
$$
|\th_j-\th_k| \ge
\begin{cases}
c_\a j^{-\a} & \text{if $j<k/2$}\\
c_\a|j-k| k^{-\a-1} & \text{if $k/2\le j\le 2k$}\\
c_\a k^{-\a} & \text{if $j>2k$}
\end{cases}
$$
where $c_\a$ is a positive constant.

For (ii), split the range of summation into two subsets: $\{(k,j):
j> \max(p,2k)\}$ and $\{(k,j): p/2<k\le p<j\le 2k\}$. The first
subset contributes at most
$$
\sumnl_{k\le p} k^{-\a-2\b}\sumnl_{j> \max(p,2k)}j^{-\a}(c_\a
k^{-\a})^{-2} = O_\ff(p^{1-\a})
$$
because $\a-2\b<-3$. The second subset contributes at most
$$
\sumnl_{p/2<k\le
p}k^{-\a-2\b}c_\a^{-2}k^{2\a+2}\sumnl_{j>p}j^{-\a}(j-k)^{-2}
=O_\ff\left(p.p^{2+\a-2\b}p^{-\a}O(1)\right),
$$
which is of order~$o_\ff(p^{-\a})$.
\end{proof}
Now remember that
$$
\ZZ_i(t) -\ZZbar(t) = \sumnl_{k\in\NN}(\tzik-\tzbar{k})\tphi_k(t)
$$
so that \bAlign \tth_k\{j=k\} &= \iint
\tK(s,t)\tphi_j(s)\tphi_k(t)\,ds\,dt
\\
&=(n-1)^{-1}\sumnl_{i\le
n}(\tz_{i,j}-\tzbar{j})(\tz_{i,k}-\tzbar{k}) , \eAlign which implies
$(n-1)^{-1}\sumnl_{i\le n} \tz_i\tz_i' = \tD^2$ and \beq
(n-1)^{-1}\sumnl_{i\le n} \teta_i\teta_i' = D^{-1}\tD^2D^{-1} :=
\diag(1,\tth_1/\th_1,\dots,\tth_N/\th_N)  . \label{teta.teta} \eeq

We will analyze~$\tK$ by rewriting it using the eigenfunctions
for~$K$. Remember that $\zij=\inner{\ZZ_i}{\phi_j}$ and the
standardized variables $\etaij=\zij/\sqrt{\th_j}$ are
independent~$N(0,1)$'s. Define~$\zbar{j}=\inner{\ZZbar}{\phi_j}$ and
$\etabar{j}=n^{-1}\sumnl_{i\le n}\etaij$ and
$$
\scov_{j,k} :=(n-1)^{-1} \sumnl_{i\le
n}\left(\etaij-\etabar{j}\right)\left(\etaik-\etabar{k}\right) ,
$$
the $(j,k)$-element of a sample covariance matrix of
i.i.d.~$N(0,I_N)$ random vectors. Then
$$
\ZZ_i(t)-\ZZbar(t) = \sumnl_{j\in\NN} (\zij-\zbar{j})\phi_j(t) =
\sumnl_{j\in\NN} \sqrt{\th_j}(\eta_{i,j}-\etabar{j})\phi_j(t)
$$
and \beq \tK(s,t) = \sumnl_{j,k\in\NN}\tK_{j,k}\phi_j(s)\phi_k(t)
\qt{with $\tK_{j,k} = \sqrt{\th_j\th_k}\scov_{j,k}$}
\label{Khat.rep} \eeq Moreover, as shown in Lemma 14 in the
supplemental Appendix, the main contribution to
${f_k=\sig_k\tphi_k-\phi_k}$ is
$$
\Lam_k := \sumnl_{j\in\NN}\Lam_{k,j}\phi_j \qt{with } \Lam_{k,j} :=
\begin{cases}
\sqrt{\th_j\th_k}\scov_{j,k}/(\th_k-\th_j)&\text{if }j\ne k\\
0&\text{if }j=k
\end{cases} .
$$
Define $$\epsilon_k := \min\{|\theta_j - \theta_k|:j\neq k\}.$$ The
following two lemmas related to perturbation theory for self-adjoint
compact operators \citep[cf. e.g.][]{Bosq2000, BirmanSolomjak87,
Kato95} are crucial in the development of Lemma \ref{txxn}. They are
special cases of Lemma 13 and Lemma 15 in the Appendix under the
general perturbation-theoretic framework. For Lemma \ref{eigenvec},
similar results were established by other authors see e.g.
\citealp[equation~2.8]{HallHosseini2006} and \citealp[Section
5.6]{CaiHall2006}. Lemma \ref{eigenspace} extends the perturbation
result for eigenprojections, obtained by \citet[Lemma 4.1]{Tyler81},
from the matrix case to the general operator case.

\begin{lemma}\label{eigenvec}
If $\epsilon_k>5\norm{\Delta}$, it follows that
$$\norm{f_k}\leq 3\norm{\Lam_k}.$$
\end{lemma}

Define $H_J = \SPAN\{\phi_j:j\in J\}$ and $\tH_J =
\SPAN\{\tphi_j:j\in J\}$ for $J\subseteq\NN$.

\Lemma\label{eigenspace} If
 $\min_{k\in J}\eps_k>5\norm{\Delta}$, then
 $$(\tH_J-H_J)\BB = \sumnl_{j\in
J}\sumnl_{k\in J^c}\phi_j b_k (\Lam_{j,k} + \Lam_{k,j}) + e$$ where
$||e||^2$ is bounded by a universal constant times $R_1 +
||\Delta||^2 R_2$ with \bAlign R_1 &= \left(\sumnl_{k\in
J}\norm{\Lam_k}^2\right) \sumnl_{k\in
J}\left(\sumnl_{j}^*\Lam_{k,j}b_j\right)^2 \\
R_2&=\sumnl_{k\in
J}\norm{\Lam_k}^2\left(\sumnl_{j}^*\frac{|b_j|}{|\th_k-\th_j|}\right)^2
+\left(\sumnl_{k\in
J}\norm{\Lam_k}|b_k|\sumnl_{j}^*\frac{1}{|\th_k-\th_j|}\right)^2\\
&\quad+\sumnl_{k\in J}\norm{\Lam_k}^2|b_k|^2k^{2+2\alpha}\eAlign
\endLemma

In fact, most of the inequalities that we need for analyzing the
estimator~$\BBhat$ defined in \cref{bhat} - \cref{N} come from
simple moment bounds (Lemma~\ref{Sjk.facts}) for the sample
covariances~$\scov_{j,k}$ and the derived bounds
(Lemma~\ref{Lam.bounds}) for the~$\Lam_k$'s.

The distribution of~$\scov_{j,k}$ does not depend on the parameters
of our model. Indeed, by the usual rotation of axes we can rewrite
$(n-1)\scov_{j,k}$ as $U_j'U_k$, where $U_1,U_2,\dots$ are
independent~$N(0,I_{n-1})$ random vectors. This representation gives
some useful  equalities and bounds. \bnote  \Lemma \label{Sjk.facts}
Uniformly over distinct $j,k,\ell$,
\begin{enumerate}
\item

$\PP \scov_{j,j}= 1$ and $\PP \left(\scov_{j,j}-1\right)^2 =
2(n-1)^{-1}$

\item
$\PP \scov_{j,k} = \PP \scov_{j,k} \scov_{j,\ell}=0 $

\item
$\PP \scov_{j,k}^2 =O(n^{-1})$

%
\end{enumerate}
\endLemma

\begin{proof} Assertion~(i) is classical because $|U_j|^2\sim\chi_{n-1}^2$.
For assertion~(ii) use $ \PP (U_1'U_2\mid U_2) =0 $ and
$$
\PP(U_1'U_2U_2'U_3\mid U_2) = \trace\left(U_2U_2'\PP
(U_3U_1')\right)=0 .
$$
For (iii) use $\PP(U_1U_1')=I_{n-1}$ and
$$
\PP(U_1'U_2U_2'U_1\mid U_2) = \trace\left(U_2U_2'\PP
(U_1U_1')\right)= \trace(U_2U_2') =|U_2|^2.
$$
\end{proof}

\Lemma  \label{Lam.bounds} Uniformly over distinct $j,k,\ell$,
\begin{enumerate}

\item\label{Lam.bounds.0}
$\PP \Lam_{k,j} = \PP \Lam_{k,j} \Lam_{k,\ell}=0 $

\item
$\PP \Lam_{k,j}^2 = O_\ff\left(n^{-1} k^{-\a}j^{-\a}
(\th_k-\th_j)^{-2}\right)$


\item \label{Lam.bounds.Lam2}
 $\PP \norm{\Lam_k}^2 =O_\ff(n^{-1}k^2)$

\end{enumerate}
\endLemma
\begin{proof} Assertions~(i) and (ii) follow from Assertions~(ii)
and~(iii) of  Lemma~\ref{Sjk.facts}. For~(iii), note that
$$
\PP \norm{\Lam_k}^2 = \sumnl_j^*\PP \Lam_{j,k}^2 =
O_\ff(n^{-1}k^{-\a})\kappa_k(2,\a)
$$
\end{proof}
\enote

To prove Lemma~\ref{txxn} we define~$\txx_{\eps,n}$ as an
intersection of sets chosen to make the six assertions of the Lemma
hold,
$$
\txx_{\eps,n} := \txx_{\Del,n}\cap\txx_{\ZZ,n}\cap \txx_{\Lam,n}\cap
\txx_{\eta,n}\cap \txx_{A,n},
$$
where the complement of each of the five sets appearing on the
right-hand side has probability less than~$\eps/5$. More
specifically, for a large enough constant~$C_\eps$, we define
\bAlign \txx_{\Del,n} &=\{\norm{\Del}\le C_\eps n^{-1/2}\}
\\
\txx_{\ZZ,n}  &=\{\maxnl_{i\le n}\norm{\ZZ_i}^2 \le C_\eps\log
n\text{ and } \,\norm{\ZZbar}\le C_\eps n^{-1/2}\}
\\
\txx_{\eta,n} &=\{\maxnl_{i\le n}|\eta_i|^2 \le C_\eps N\log n\}
\qt{as in Section~\ref{known}}
\\
\txx_{A,n}&= \{\twonorm{\sumnl_{i\le n} \teta_i\teta_i' } \le C_\eps
n\} \eAlign The definition of~$\txx_{\Lam,n}$, in
subsection~\ref{proof.H}, is slightly more complicated.  It is
defined by requiring various functions of the~$\Lam_k$'s to be
smaller than~$C_\eps$ times their expected values.

The set~$\txx_{A,n}$ is almost redundant. From
Definition~\ref{ff.def} we know that
$$
\min_{1\le j<j'\le N}|\th_j-\th_{j'}| \ge (\a/R) N^{-1-\a}
\AND/\min_{1\le j\le N}\th_j \ge R^{-1}N^{-\a}  .
$$
The choice~$N\sim n^\zeta$ with $\zeta< (2+2\a)^{-1}$ ensures
that~$n^{1/2}N^{-1-\a}\to\infty$. On~$\txx_{\Del,n}$ the spacing
assumption used in Lemmas \ref{eigenvec} and \ref{eigenspace} holds
for all~$n$ large enough; all the bounds from those lemmas are
available to us on~$\txx_{\eps,n}$. In particular,
$$
\maxnl_{j\le N}|\tth_j/\th_j-1| \le O_\ff(N^\a\norm{\Del}) =o_\ff(1)
.
$$
Equality~\cref{teta.teta} shows that $\txx_{A,n}\subseteq
\txx_{\Del,n}$ eventually if we make sure $C_\eps>1$.

\subsubsection{Proof of Lemma~\ref{txxn} part \ref{txxn.Del}}\label{txxnDel} Observe
that \bAlign \PP \norm{\Del}^2 &=
\sumnl_{j,k}\PP\left(\tK_{j,k}-\th_j\{j=k\}\right)^2
= \sumnl_{j,k}\th_j\th_k\PP\left(S_{j,k}-\{j=k\}\right)^2\\
&\le \sumnl_{j}\th_jO_\ff(n^{-1}) + \sumnl_{j,k}\th_j\th_k
O_\ff(n^{-2})= O_\ff(n^{-1}) \eAlign

\subsubsection{Proof of Lemma~\ref{txxn} part \ref{txxn.ZZ}}\label{txxnZZ} As
before, Lemma~\ref{gaussian} controls~$\maxnl_{i\le
n}\norm{\ZZ_i}^2$. To control the $\ZZbar$ contribution, note that
$n\norm{\ZZbar}^2$ has the same distribution as $\norm{\ZZ_1}^2$,
which has expected value $\sumnl_{j\in\NN}\th_j<\infty$.

\subsubsection{Proof of Lemma~\ref{txxn} parts \ref{txxn.Hm} and
\ref{txxn.HN}}\label{txxnH} \label{proof.H} Calculate expected
values for all the terms that appear in the bound of Lemma
\ref{eigenspace}. \bAlignL \PPnmuK&\sumnl_{k\le
p}\left(\sumnl_{j>p}\Lam_{k,j}b_j\right)^2 +
\PPnmuK\sumnl_{j>p}\left(\sumnl_{k\le p}\Lam_{k,j}b_k\right)^2
\notag\\
&= \sumnl_{k\le
p}\sumnl_{j>p}\PPnmuK\Lam_{k,j}^2\left(b_j^2+b_k^2\right) \qt{by
Lemma~\ref{Lam.bounds} part \ref{Lam.bounds.0}}
\notag\\
&= O_\ff(n^{-1})\sumnl_{k\le
p}\sumnl_{j>p}k^{-\a-2\b}j^{-\a}(\th_k-\th_j)^{-2}
\notag\\
&= O_\ff(n^{-1}p^{1-\a}) \qt{by Lemma~\ref{weights}} \label{proj1}
\eAlignL and
\bnote \bAlign ||\Delta||^2\PPnmuK \sumnl_{k\le
p}b_k^2\norm{\Lam_k}^2 k^{2+2\alpha} &=
O_\ff(n^{-1}||\Delta||^2)\sumnl_{k\le p}k^{4+2\alpha-2\b}\\& =
O_\ff(n^{-2})\left(1+p^{5+2\alpha-2\b}+\log p\right) \eAlign \enote
and
$$
\PPnmuK \sumnl_{k\le p}|b_k|\norm{\Lam_k}^2 =
O_\ff(n^{-1})\sumnl_{k\in J}k^{2-\b} =
O_\ff(n^{-1})\left(1+p^{3-\b}+\log p\right)
$$
and
$$
\PPnmuK \sumnl_{k\le p}\norm{\Lam_k}^2 = O_\ff(n^{-1}p^3)
$$
and \bAlignL \PPnmuK \sumnl_{k\le
p}\left(\sumnl_{j}^*\Lam_{k,j}b_j\right)^2 &=
O_\ff(n^{-1})\sumnl_{k\le p}\sumnl_{j}^*
k^{-\a}j^{-a-2\b}(\th_k-\th_j)^{-2}
\notag\\
&= O_\ff(n^{-1}) \qt{by Lemma~\ref{weights}} \label{proj2} \eAlignL
and \bAlignL ||\Delta||^2\PPnmuK&\sumnl_{k\le p}
\norm{\Lam_k}^2\left(\sumnl_j^*\frac{|b_j|}{|\th_k-\th_j|}\right)^2
\notag\\
&= O_\ff(n^{-1}||\Delta||^2)\left(p^3+ p^{5+2\a-2\b}\log^2 p\right)
\label{proj3} \eAlignL and \beq \sumnl_{k\le
p}b_k^2\left(\sumnl_j^*\frac{1}{|\th_k-\th_j|}\right)^2
=O_\ff(1+p^{3+2\a-2\b}\log^2 p) \qt{by Lemma~\ref{weights}}.
\label{proj4} \eeq For some constant $C_\eps=C_\eps(\ff)$, on a
set~$\xx_{\Lam,n}$ with $\PPnmuK \xx_{\Lam,n}^c <\eps$, each of the
random quantities in the previous set of inequalities (for both
$p=m$ and $p=N$) is bounded by $C_\eps$ times its $\PPnmuK$ expected
value. By virtue of Lemma~\ref{Lam.bounds} part
\ref{Lam.bounds.Lam2}, we may also assume that~$\norm{\Lam_k}^2\le
C_\eps k^2/n$ on~$\xx_{\Lam,n}$.

\bnote From Lemma \ref{eigenspace}, it follows that on the set
$\xx_{\Del,n}\cap \xx_{\Lam,n}$, if $p\le N$, \bAlign
&\norm{(\tH_p-H_p)\BB}^2\\
&\qquad\le O_\ff(n^{-1}p^{1-\a})
+ O_\ff(n^{-2}) \left(1+p^{5+2\alpha-2\b}+\log p +p^{6-\b}+\log^2p\right)\\
&\qquad\qquad+ O_\ff(n^{-1}p^3) O_\ff(n^{-1}) +
O_\ff(n^{-2})\left(p^3+ p^{5+2\a-2\b}\log^2 p\right)
\\
&\qquad\qquad+ O_\ff(n^{-2}p^3)O_\ff(1+p^{3+2\a-2\b}\log^2 p)
 \\
 &\qquad= O_\ff(n^{-1}p^{1-\a})
\eAlign This inequality leads to the asserted conclusions when~$p=m$
or~$p=N$. \enote

\subsubsection{Proof of Lemma~\ref{txxn} part \ref{txxn.teta}}\label{txxnteta} By
construction, $\teta_{i1}=1$ for every~$i$ and, for $j\ge 2$,
$$
\sqrt{\theta_j}\teta_{i,j} =(\tz_{i,j}-\tzbar{j}) =
\inner{\ZZ_i-\ZZbar}{\tphi_j}
$$
Thus, for $j\ge2$,
$$
\sig_j\teta_{i,j} =
\th_j^{-1/2}\inner{\ZZ_i-\ZZbar}{\phi_j+f_j}=\eta_{i,j}+\tdel_{i,j}
$$
with, due to Lemma \ref{eigenvec},
$$
|\del_{i,j}|^2\le
\th_j^{-1}\left(\norm{\ZZ_i}+\norm{\ZZbar}\right)^2\norm{f_j}^2 \le
O_\ff\left(\frac{j^{2+\a}\log n}{n}\right) \qt{on $\txx_{\eps,n}$}.
$$
In vector form, \beq \tS\teta_i = \eta_i + \tdel_i \qt{with
}|\tdel_i|^2 =O_\ff\left(\frac{N^{3+\a}\log n}{n}\right)\le
o_\ff(n/N^2) \text{ on }\txx_{\eps,n}. \label{teta.eta} \eeq It
follows that
$$
\maxnl_{i\le n}|\teta_i| = \maxnl_{i\le n}|\tS\teta_i| \le
\maxnl_{i\le n}|\eta_i|+o_\ff(\sqrt{n}/N) = O_\ff(\sqrt{n}/N) \qt{on
$\txx_{\eps,n}$}.
$$

\subsubsection{Proof of Lemma~\ref{txxn} part \ref{txxn.tA}}\label{txxntA} From inequality~\cref{tlam.lam} we know that
$$
\tilde{\eps}_{N} := \maxnl_{i\le n}|\tlam_{i,N}-\lam_{i,N}| =
O_\ff(n^{-(1+\nu')/2}) \qt{on $\txx_{\eps,n}$}
$$
and from the Section \ref{known} we have $ \maxnl_{i\le
n}|\lam_{i,N}| = O_\ff(\rootlogn) $. Assumption~(\psiPPP) in
Section~\ref{model} and the Mean-Value theorem then give
$$
\maxnl_{i\le n}|\psiddot(\tlam_{i,N})-\psiddot(\lam_{i,N})| \le
\tilde{\eps}_N\psiddot(\lam_{i,N})G(\tilde{\eps}_N) =o_\ff(1) .
$$
If we replace $\psiddot(\tlam_{i,N})$ in the definition of~$\tA_n$
by~$\wi:=\psiddot(\lam_{i,N})$ we make a change~$\Gamma$ with
$$
\twonorm{\Gamma} \le o_\ff(1)\twonorm{(n-1)^{-1}\sumnl_{i\le n}
\teta_i\teta_i' },
$$
which, by equality~\cref{teta.teta}, is of order~$o_\ff(1)$
on~$\txx_{\eps,n}$.

From Assumption~(\psiPP) we have $c_n :=\log\maxnl_{i\le n}\wi=
o_\ff(\log n)$. Uniformly over all unit vectors~$u$ in~$\RR^{N+1}$
we therefore have \bAlign u'\tS\tA_n\tS u &= o_\ff(1) +
(n-1)^{-1}\sumnl_{i\le n}\wi u'(\eta_i+\tdel_i)(\eta_i+\tdel_i)'u
\\
&= o_\ff(1) + \left(1+O(n^{-1})\right)u'A_nu
\\
&\qquad+O_\ff\left(n^{-1}\right)\sumnl_{i\le n}\wi
\left((u'\tdel_i)^2 +
 2(u'\eta_i)(u'\tdel_i)\right)
\eAlign Rearrange then take a supremum over~$u$ to conclude that
$$
\twonorm{\tS\tA_n\tS-A_n} \le o_\ff(1) + O_\ff(e^{c_n})\maxnl_{i\le
n}\left(|\tdel_i|^2 + 2|\tdel_i|\,|\eta_i|\right)
$$
Representation~\cref{teta.eta} and the defining property
of~$\txx_{\eta,n}$ then ensure that the upper bound is of
order~$o_\ff(1)$ on~$\txx_{\eps,n}$.

\bibliographystyle{imsart-nameyear}
\bibliography{DPZ_revision_arXiv}

\newpage

\section{Appendix}\label{appendix}In this section, we introduce some useful results in spectral theory and perturbation theory.
Some of the results are well-established. We briefly review them for
the purpose of easy reference. For example, the results for
eigenvalues have become quite standard for decades \citep[see,
e.g.][Chapter VII.6]{DunfordSchwartz88}. We derive a bound for the
perturbation of eigenprojections (Lemma 15) which plays a key role
in the slope function estimation problem. This bound is closely
related Proposition 2 in \citet[]{CardotMasSarda2007}, which was
tailored to solve the prediction problem at a random design.
However, the two results are different. A comparison between their
result and our bound in Lemma~\ref{HBB} is discussed later following
Lemma~\ref{HBB}. We could not find the same (or stronger) bound
explicitly in the existing perturbation literature.

The spectral theory and the perturbation theory in Hilbert spaces
have been serving as powerful tools that allow statisticians to
tackle the statistical approximation problems in an elegant way.
From Lemma \ref{spectral.values} to Lemma \ref{spectral.vectors.2}
we shall review the well-established perturbation-theoretic results
for eigenvalues and eigenvectors of positive and self-adjoint
compact operators respectively. Our main contribution of this
section is to extend the perturbation result for eigenprojections,
obtained by \citet[Lemma 4.1]{Tyler81}, from the matrix case to the
general operator case. Our perturbation result for eigenprojections
will be introduced in Lemma \ref{HBB}.

Suppose $T$ is a positive and self-adjoint compact operator in a
Hilbert space $\hh$. According to the spectral theory for positive
and self-adjoint compact operators \citep[see e.g.][Page
209]{BirmanSolomjak87}, the operator $T$ has a sequence of
decreasing nonnegative eigenvalues $\{\theta_i\}$ and a sequence of
corresponding eigenvectors $\{e_i\}$. That is, $Te_i = \theta_i e_i$
with $\theta_1 \geq \theta_2 \geq \cdots \geq 0$. Furthermore, $T$
has the spectral decomposition \beq\label{spectral.decomp} T =
\sumnl_{k\in \NN} \theta_k e_k \otimes e_k \eeq which converges in
the operator norm.


In this section, the perturbation-theoretic results are the
functional analysis results without involving randomness. The focus
here is on the results for positive and self-adjoint compact
operators, and for more general discussion on perturbation theory of
linear operators please see, for example, \citet[]{Kato95}. More
precisely, let $\tT$ be another positive and self-adjoint compact
operator in $\hh$ with spectral decomposition
\beq\label{tT.spectral.decomp} \tT = \sumnl_{k\in \NN}\tth_k \te_k
\otimes \te_k \eeq The eigenprojection of the operator $\tT$
associated with eigenvalues $\widetilde{\Theta}_J := \{\tth_j:j\in
J\}$, denoted by $\tH_J$, is the orthogonal projection onto the
eigenspace of $\tT$ associated with $\widetilde{\Theta}_J$, that is,
$\SPAN\{\te_j~:~j\in J\}$. In fact, we have $\tH_J = \sumnl_{j\in
J}\te_j \otimes \te_j$. Analogously, the eigenprojection $H_J$ can
be defined for the operator $T$. We shall study how well the
differences of $\theta_j - \tth_j$, $e_j - \te_j$, and $H_J - \tH_J$
can be controlled by $\Delta = T - \tT$, given that $\delta :=
||\Delta||_2$ is small.

In statistical applications, the operator $T$ is usually taken as
unknown, while the operator $\tT$ taken as the estimation of $T$.
Perturbation theory suggests that as long as $\tT$ approximates $T$
well, the eigen-elements of $\tT$ can project the analogous
eigen-elements of $T$ well. This idea has been explored and utilized
by \citet[]{Tyler81}, \citet[]{Bosq2000}, \citet[]{CaiHall2006},
\citet[]{HallHorowitz2007}, and \citet[]{CardotMasSarda2007}, among
others. More interestingly, \citet[]{HallHosseini2006} proposes a
Taylor-expansion type of approximation of eigenvectors which is
better adapted to the statistical approximation purposes.

In the application of Section \ref{unknown}, we draw probabilistic
conclusions when~$\tT$ is random for the special case where $T=\KK$,
the population covariance kernel, and $\tT=\tKK$, the sample
covariance kernel, both acting on $\hh=\ll^2[0,1]$. The eigenvectors
$\{e_i\}$ and $\{\te_i\}$ will be principal components $\{\phi_i\}$
and $\{\tphi_i\}$ respectively.

Before formally illustrating the perturbation-theoretic results in
details, we shall introduce some necessary notations and basic
mathematical relations here. Because $\{e_j: j\in \NN\}$ forms a
complete orthonormal basis for the Hilbert space $\hh$, the operator
$\tT$ also has the following representation \beq\label{tT} \tT =
\sumnl_{j,k\in \NN}\tT_{j,k} e_j \otimes e_k \eeq which converges in
the operator norm.

Note that $\tT_{j,k}=\tT_{k,j}$ because~$\tT$ is self-adjoint. This
representation gives
$$
\Del =
\sumnl_{j,k\in\NN}\left(\tT_{j,k}-\th_j\{j=k\}\right)e_j\otimes e_k
$$
and
$$
\delta^2=\norm{\Del}^2 = \sup\NL_{\norm{x}=1}\inner{x}{\Del x}^2 \le
\sumnl_{j,k\in\NN}\left(\tT_{j,k}-\th_j\{j=k\}\right)^2 .
$$

We also define \beq\label{Lambda} \Lam_k :=
\sum_{j\in\NN}\Lam_{k,j}e_j \qt{with } \Lam_{k,j} :=
\begin{cases}
\tT_{j,k}/(\th_k-\th_j)&\text{if }j\ne k\\
0&\text{if }j=k
\end{cases} .
\eeq

Notice that $\{\te_k:k\in\NN\}$ is also an orthonormal basis
for~$\hh$. Define ${\sig_{j,k}:=\inner {e_j} {\te_k}}$. Then
$$
e_j = \sumnl_{k\in\NN}\sig_{j,k}\te_k \AND/ \te_k =
\sumnl_{j\in\NN}\sig_{j,k}e_j
$$
and
$$
\{j = j'\} = \inner{e_j}{e_{j'}} = \sumnl_{k\in
\NN}\sig_{j,k}\sig_{j',k} .
$$

We cannot hope to find a useful bound on~$\norm{\te_k-e_k}$, because
there is no way to decide which of $\pm\te_k$ should be
approximating~$e_k$. However, we can bound $\norm{f_k}$, where
\beq\label{fk} f_k = \sig_k\te_k -e_k \qt{with } \sig_k :=
\sign\left(\sig_{k,k}\right) :=
\begin{cases}
+1&\text{if }\sig_{k,k}\ge0\\
-1&\text{otherwise}
\end{cases},
\eeq which will be enough for our purposes.

To simplify notations, write $\sumnl_j^*$ for
$\sumnl_{j\in\NN}\{j\ne k\}$ and $\sumnl_k^*$ for
$\sumnl_{k\in\NN}\{k\ne j\}$ in this section.

The following lemma has been proved in multiple places
. Here we stated it with a brief
proof for the purpose of easy reference.
\begin{lemma}\label{spectral.values}
Suppose $T$ and $\tT$ are two positive and self-adjoint operators
with spectral decompositions \cref{spectral.decomp} and
\cref{tT.spectral.decomp}, then it follows that
 \beq |\th_j-\tth_j| \le
\del\qt{for all $j\in\NN$}. \label{evalue.approx} \eeq
\end{lemma}
\begin{proof}
The eigenvalues have a variational characterization; see
\citet[Section 4.2]{Bosq2000} or \citet[Chapter
9]{BirmanSolomjak87}: \beq \th_j = \inf_{\dim(L)<j} \sup\{\inner x
{Tx}: x\perp L\text{ and } \norm{x}=1 \} . \label{var.evalue} \eeq
The first infimum runs over all subspaces~$L$ with dimension at
most~$j-1$. (When~$j$ equals~$1$ the only such subspace
is~$\emptyset$.) Both the infimum and the supremum are achieved: by
$L_{j-1}=\SPAN\{e_i: 1\le i <j\}$ and $x=e_j$. Similar assertions
hold for~$\tT$ and its eigenvalues.

By the analog of~\cref{var.evalue} for~$\tT$,  \bAlign
\tth_j &\leq \sup\{\inner x {\tT x}: x\perp L_{j-1}\text{ and } \norm{x}=1 \} \\
&\leq \sup\{\inner x {T x} +\del: x\perp L_{j-1}\text{ and }
\norm{x}=1 \} = \th_j+\del . \eAlign   Argue similarly with the
roles of~$T$ and~$\tT$ reversed to conclude the result.
\end{proof}

In order to approximate an eigenvector $e_k$ reasonably well, we
need to assume that the eigenvalue~$\th_k$ is well separated from
the other~$\th_j$'s, to avoid the problem that the eigenspace
of~$\tT$ for the eigenvalue~$\tth_k$ might have dimension greater
than one. More precisely, we consider a~$k$ for which
\beq\label{def.epsilon.k} \eps_k := \min\{|\th_{j}-\th_k|: j\ne k\}
>5\delta ,
\eeq which implies
$$
|\tth_k-\th_j|\ge |\th_k-\th_j| -\delta \ge \tfrac45|\th_k-\th_j|
\ge \tfrac45\eps_k .
$$

The following lemmas provide approximative results for $f_k$ under
the assumption that ${\epsilon_k
> 5\delta}$. Similar results were established by other authors; see,
for example, \citealp[equation~2.8]{HallHosseini2006} and
\citealp[Section 5.6]{CaiHall2006}.
\begin{lemma}\label{spectral.vectors.1}
Suppose $T$ and $\tT$ are two positive and self-adjoint operators
with spectral decompositions \cref{spectral.decomp} and
\cref{tT.spectral.decomp}. The vectors $\{\Lam_k\}$ and $\{f_k\}$
are defined as in \cref{Lambda} and \cref{fk} respectively. Then if
$\epsilon_k>5\delta$, it follows that
$$\norm{f_k}\leq 3\norm{\Lam_k}.$$
\end{lemma}
\begin{proof}
The starting point for our approximations is the equality \beq
\inner{\Del \te_k}{e_j} = \inner{\tT\te_k}{e_j} -
\inner{\te_k}{Te_j} = (\tth_k-\th_j)\sig_{j,k} . \label{Dee} \eeq
For $j\ne k$ we then have
$$
\frac{16}{25}(\th_k-\th_j)^2\sig_{j,k}^2 \le \inner{\sig_k\Del
\te_k}{e_j}^2
 \le 2 \inner{\Del f_k}{e_j}^2 + 2 \inner{\Del e_k}{e_j}^2 ,
$$
which implies
$$
\sig_{j,k}^2 \le \frac{25}{8}\inner{\Del f_k}{e_j}^2/\eps_k^2  +
2\tT_{j,k}^2/(\th_k-\th_j)^2 \qt{because $\inner{T e_k}{e_j}=0$ for
$j\ne k$.}
$$

The introduction of the $\sig_k$ also ensures that \bAlign
\norm{f_k}^2 &= \norm{e_k}^2+\norm{\te_k}^2 -
2\sig_k\inner{e_k}{\te_k} = 2-2|\sig_{k,k}|
\\
&\le 2-2\sig_{k,k}^2 \qt{because $|\sig_{k,k}|\le 1$}
\\
& =  2\sumnl_j^*\sig_{j,k}^2
\\
&\le\sumnl_j^*\frac{25}{4}\inner{\Del f_k}{e_j}^2/\eps_k^2
+\frac{25}{4}\sumnl_j^*\tT_{j,k}^2/(\th_k-\th_j)^2 . \eAlign The
first sum on the right-hand side is less than
$$
\frac{25}{4} \norm{\Del f_k}^2 /\eps_k^2 \le \delta^2\norm{f_k}^2
/(4\del^2) =  \norm{f_k}^2 /4.
$$
The second sum can be written as $25\norm{\Lam_k}^2/4$. Then,

\beqN \norm{f_k}^2\leq \frac{25}{3}\norm{\Lam_k}^2 <
9\norm{\Lam_k}^2
.\eeqN
\end{proof}
\begin{lemma}\label{spectral.vectors.2} Suppose $T$ and $\tT$ are two positive and self-adjoint operators
with spectral decompositions \cref{spectral.decomp} and
\cref{tT.spectral.decomp}. The vectors $\{\Lam_k\}$ and $\{f_k\}$
are defined as in \cref{Lambda} and \cref{fk} respectively. Then if
$\epsilon_k>5\delta$, the corresponding operator $f_k$ has the
representation:
$$f_k = \Lam_k + r_k$$
with
$$\inner{r_k}{e_k}=-\frac{1}{2}\norm{f_k}^2~~~\mbox{and}~~~|\inner{r_k}{e_j}|\leq
\frac{5\delta \norm{\Lam_k}}{|\theta_k - \theta_j|}~~~\forall~~j\neq
k.$$
\end{lemma}
\begin{proof}
Start once more from equality~\cref{Dee}. For $j\ne k$, \bAlignL
\sig_k\sig_{j,k} &= \sig_{k} \inner{\Del \te_k}{e_j}/(\tth_k-\th_j)
\notag\\
&= \inner{\Del (e_k+f_k)}{e_j}/(\th_k + \gam_k-\th_j) \qt{where
$\gam_k=\tth_k-\th_k$}
\notag\\
&= \Lam_{k,j}\left(1 -\frac{\gam_k}{\th_j-\th_k}\right)^{-1} +
\frac{ \inner{\Del f_k}{e_j} }{\tth_k-\th_j} \qt{because
$\inner{Te_k}{e_j}=0$}
\notag\\
&= \Lam_{k,j} +r_{k,j} \qt{where } r_{k,j}:=
\frac{\tth_k-\th_k}{\th_j-\th_k}\Lam_{k,j} + \frac{ \inner{\Del
f_k}{e_j} }{\tth_k-\th_j}. \label{LamR}\eAlignL The $ r_{k,j}$'s are
small: \bAlignL |r_{k,j}| &\le
\frac54\left(\frac{\del|\Lam_{k,j}|+|\inner{\Del
f_k}{e_j}|}{|\th_k-\th_j|}\right) \qt{for $j\ne k$, if
$\eps_k>5\del$}
\notag\\
&\le \frac{5\del\norm{\Lam_k}}{|\th_k-\th_j|} \qt{by Lemma
\ref{spectral.vectors.1}}. \label{rjk} \eAlignL Define
${r_{k,k}=|\sig_{k,k}|-1= -\tfrac12\norm{f_k}^2}$ and $r_k=
\sumnl_{j\in\NN}r_{k,j}e_j$. Combine \cref{Lambda} and \cref{LamR},
we then have the representation: \beq f_k = \sig_k\te_k-e_k =
\left(\sig_k\inner{\te_k}{e_k} -1\right)e_k
+\sumnl_j^*\sig_k\sig_{j,k}e_j = \Lam_k + r_k  . \label{fk.rep} \eeq
\end{proof}

In the rest of this section, we shall establish an approximation for
$\tH_J\BB-H_J\BB$ for a $\BB=\sumnl_{j} b_je_j$ in~$\hh$, an
extension of the finite-dimensional perturbation result
\citet[][Lemma 4.1]{Tyler81} to the case of general
infinite-dimensional operators.

The difference $\tH_J-H_J$ equals \bAlignL \sumnl_{k\in
J}&(\sig_k\te_k)\otimes (\sig_k\te_k)-e_k\otimes e_k
\notag\\
&= \sumnl_{k\in J}\sig_k\te_k\otimes r_k + \sumnl_{k\in
J}(e_k+f_k)\otimes \Lam_k
\notag\\
&\quad +\sumnl_{k\in J}\left((e_k+\Lam_k +r_k)\otimes e_k -
e_k\otimes e_k\right)
\notag\\
&= \rr_J +\sumnl_{k\in J} e_k\otimes\Lam_k +\Lam_k\otimes e_k \notag\\
&\qt{where } \rr_J := \sumnl_{k\in J} \sig_k\te_k\otimes r_k +
f_k\otimes \Lam_k + r_k\otimes e_k .\label{rr} \eAlignL

Self-adjointness of~$\tT$  implies  $\tT_{j,k}=\tT_{k,j}$ and hence
$\Lam_{j,k}=-\Lam_{k,j}$. The anti-symmetry eliminates some terms
from the main contribution to~$\tH_J-H_J$: \bAlignL \sumnl_{k\in J}
&e_k\otimes\Lam_k +\Lam_k\otimes e_k = \sumnl_{k\in J}\sumnl_{j\in
J^c} \Lam_{k,j}\left(e_k\otimes e_j +e_j\otimes e_k\right) .
\eAlignL With this simplification we get the following
representation for~$(\tH_J-H_J)\BB$: \beqN (\tH_J-H_J)\BB =
\sumnl_{j\in J}\sumnl_{k\in J^c}e_j b_k (\Lam_{j,k} + \Lam_{k,j}) +
\rr_J\BB. \eeqN
\bnote For the three contributions to the bound
for~$\norm{\rr_J\BB}^2$ we make repeated use of the inequalities,
based on Lemma \ref{spectral.vectors.1} and Lemma
\ref{spectral.vectors.2}, \bAlignL |\inner {r_k}{x}| &\le
\frac{3}{2}\norm{\Lam_k}\norm{f_k}|x_k| +
5\del\norm{\Lam_k}\sumnl_j^* \frac{|x_j|}{|\th_k-\th_j|}\label{rx1}\\
&\leq \frac{9}{2}\norm{\Lam_k}^2|x_k| + 5\del\norm{\Lam_k}\sumnl_j^*
\frac{|x_j|}{|\th_k-\th_j|} \label{rx2}\eAlignL which is valid
whenever $\eps_k>5\del$. Combine \cref{rx1} and the following
well-known inequality \citep[see e.g.][Equ. 5.2]{HallHorowitz2007}:
$$||f_k||\leq 2\sqrt{2}\delta \min\{\theta_{k-1}-\theta_k, \theta_k-\theta_{k+1}\}^{-1}=2\sqrt{2}\delta\epsilon_k^{-1},$$
we get $$|\inner{r_k}{x}|\leq
C\delta\norm{\Lam_k}\left(\epsilon_k^{-1}|x_k|+\sumnl_j^*
\frac{|x_j|}{|\th_k-\th_j|}\right).$$ To avoid an unnecessary
calculation of precise constants, we adopt the convention of the
variable constant: we write~$C$ for a universal constant whose value
might change from one line to the next. The first two contributions
are: \bAlign \norm{\sumnl_{k\in J}\sig_k\te_k \inner{r_k}{\BB} }^2
&= \sumnl_{k\in J}\inner {r_k}{\BB}^2
\\
&\le C \delta^2\sumnl_{k\in J}\norm{\Lam_k}^2 \left[ b_k^2
\epsilon_k^{-2}
+\left(\sumnl_j^*\frac{|b_j|}{|\th_k-\th_j|}\right)^2\right] \eAlign
and \bAlign \norm{\sumnl_{k\in J}f_k\inner{\Lam_k}{\BB}}^2 &\le
\left(\sumnl_{k\in J}\norm{f_k}\,|\inner{\Lam_k}{\BB}|\right)^2
\\
&\le C\left(\sumnl_{k\in J}\norm{\Lam_k}^2\right) \sumnl_{k\in
J}\left(\sumnl_{j}^*\Lam_{k,j}b_j\right)^2. \eAlign For the third
contribution, let $x=\sumnl_j x_j e_j$ be an arbitrary unit vector
in~$\hh$. Then \bAlignL &\left(\sumnl_{k\in J}\inner{r_k\otimes e_k
\BB}{x}\right)^2 =\left(\sumnl_{k\in J}b_k\inner{r_k}{x}\right)^2
\notag\\
&\quad\le C\delta^2\left[\sumnl_{k\in
J}|b_k|\norm{\Lam_k}\left(|x_k|\epsilon_k^{-1}
+\sumnl_j^*\frac{|x_j|}{|\th_k-\th_j|}\right)\right]^2 \label{xbb}
\\
&\quad\le C\delta^2\sumnl_{k\in J}|b_k|^2\norm{\Lam_k}^2
\epsilon_k^{-2}+C\del^2 \left(\sumnl_{k\in
J}\norm{\Lam_k}|b_k|\sumnl_j^*\frac{1}{|\th_k-\th_j|}\right)^2
\label{normbound}.\eAlignL take the supremum over~$x$, which doesn't
even appear in the last line, to get the same bound for
$\norm{\sumnl_{k\in J}b_kr_k}^2$.

In sum, we can obtain the following lemma: \Lemma\label{HBB} If
 $\min_{k\in J}\eps_k>5\delta$, then
 $$(\tH_J-H_J)\BB = \sumnl_{j\in
J}\sumnl_{k\in J^c}e_j b_k (\Lam_{j,k} + \Lam_{k,j}) + \rr_J \BB$$
where $\rr_J$ is defined in \cref{rr} and $||\rr_J \BB||^2$ is
bounded by a universal constant times $R_1 + \delta^2 R_2$  with
\bAlign R_1 &= \left(\sumnl_{k\in J}\norm{\Lam_k}^2\right)
\sumnl_{k\in
J}\left(\sumnl_{j}^*\Lam_{k,j}b_j\right)^2 \\
R_2&=\sumnl_{k\in
J}\norm{\Lam_k}^2\left(\sumnl_j^*\frac{|b_j|}{|\th_k-\th_j|}\right)^2
+\left(\sumnl_{k\in
J}\norm{\Lam_k}|b_k|\sumnl_j^*\frac{1}{|\th_k-\th_j|}\right)^2\\
&\quad+\sumnl_{k\in J}\norm{\Lam_k}^2|b_k|^2\epsilon_k^{-2}\eAlign
\endLemma
\enote

This lemma is the keystone to establish the parts~(iii) and~(iv) in
Lemma~\ref{txxn}. It is similar to Proposition 2 in
\citet{CardotMasSarda2007} in the sense that both deal with the the
approximation problems of eigenprojections. Particularly, we observe
that the same trick of using anti-symmetry is applied to eliminate
some terms from the main contributions to the approximation errors
(see Equation (45) in this section and Equation (23) in
\citet{CardotMasSarda2007}). However, the two results are different.
The other authors consider the bound for ${\inner{(\tH_J -
H_J)\BB}{\XX_{n+1}}}$ which is motivated by the prediction problem
at a random design, whereas we establish a bound for ${||(\tH_J -
H_J)\BB||}$ which is relevant to the slope function estimation
problem. More precisely, the independent randomness of $\XX_{n+1}$
helps cancel off many cross-product terms and accelerates the decay
rates of the summands in the expansion. See, for example, Equation
(24) and (25) in~\citet{CardotMasSarda2007}. Due to the `smoothing'
effect of the independent random curve $\XX_{n+1}$, we cannot
directly apply the convergence result in~\citet[Proposition
2]{CardotMasSarda2007} to our case. Besides, the bound in the lemma
above is a pure mathematical perturbation-theoretic result not
involving any randomness treatment, which we believe is a
potentially more general result.

\end{document}